\documentclass[11pt]{amsart}
\usepackage{amsmath}

\numberwithin{equation}{section}

\theoremstyle{plain}

\newtheorem{MainThm}{Main Theorem}

\newtheorem{Prop}{Proposition}[section]
\newtheorem{Cor}[Prop]{Corollary}
\newtheorem{Lem}[Prop]{Lemma}
\newtheorem{Def}[Prop]{Definition}
\newtheorem{Thm}[Prop]{Theorem}
\newtheorem{Rem}[Prop]{Remark}

\theoremstyle{remark}

\newcommand{\CBbb}{\mathbb C}
\newcommand{\RBbb}{\mathbb R}
\newcommand{\ZBbb}{\mathbb Z}

\newcommand{\ad}{\mathop{\rm ad}\nolimits }
\newcommand{\rk}{\mathop{\rm rk}\nolimits }
\newcommand{\tr}{\mathop{\rm Tr}\nolimits }

\newcommand{\hol}{\mathop{\rm hol}\nolimits }

\newcommand{\met}{\mathop{\rm Met}_{-1}\nolimits }

\newcommand{\I}{\mathop{\mathbb I}\nolimits }
\newcommand{\diag}{\mathop{\rm diag}\nolimits }
\newcommand{\Hom}{\mathop{\rm Hom}\nolimits }

\newcommand{\dbar}{\bar\partial}
\newcommand{\lra}{\longrightarrow}

\newcommand{\A}{\mathcal{A}}

\newcommand{\Ab}{{\A}_\beta}
\newcommand{\R}{\mathcal{R}}
\newcommand{\Ra}{{\R}_\alpha}
\newcommand{\Rb}{{\R}_\beta}
\newcommand{\Ma}{{\M}_\alpha}

\newcommand{\T}{\mathcal{T}}
\newcommand{\M}{\mathcal{M}}

\newcommand{\G}{\mathfrak{G}}
\newcommand{\Gc}{{\G}^{\CBbb}}
\newcommand{\ab}{{\alpha\beta}}
\newcommand{\pisigma}{\pi^{[\sigma^\ast]}}

\newcommand{\foliation}{\mathcal{F}^{[\sigma^\ast]}}

\newcommand{\foliationtau}{\mathcal{F}^{[\sigma^\ast(\ell)]}}
\newcommand{\foliationphi}{\mathcal{F}^{[\sigma^\ast(0)]}}
\newcommand{\altA}{\underset{\sim}{A}}
\newcommand{\alth}{\underset{\sim}{h}}
\newcommand{\altg}{\underset{\sim}{g}}

\begin{document}

% Topmatter

\title[The Yang-Mills Flow]
	{ The Yang-Mills Flow Near the  Boundary \\ of Teichm\"uller Space}

\bigskip

\author[G. D. Daskalopoulos]{Georgios D. Daskalopoulos}

\address{Department of Mathematics \\
		Brown University \\
		Providence,  RI  02912}

\email{daskal@gauss.math.brown.edu}

\thanks{Supported in part by NSF grant DMS-9504297}

\author[R. A. Wentworth]{Richard A. Wentworth}

\address{Department of Mathematics \\
   University of California \\
   Irvine,  CA  92697}

\email{raw@math.uci.edu}

\thanks{Supported in part by NSF grant DMS-9503635 and a Sloan Fellowship}

\date{May 11, 1998}

% End Topmatter

\begin{abstract}  We study the behavior of the Yang-Mills flow for unitary
connections on compact and non-compact oriented surfaces with varying metrics. The 
flow can be used to define a one dimensional foliation on the space of $SU(2)$ representations of
a once punctured surface.  This foliation universalizes over Teichm\"uller space and is
equivariant with respect to the action of the mapping class group. It is shown how to extend the
foliation  as a singular foliation over the Strebel boundary of Teichm\"uller space, and continuity
of this extension is the main result of the paper.
\end{abstract}

\maketitle

\baselineskip=16pt

%%%%%%%%%%%%%%%%%%%%%%%%%%%%%%%%%%%%%%%%%%%%%%%%%%%%%%%%%%%%%%%%%%%%%%%%%%%%%

\section{Introduction}

%%%%%%%%%%%%%%%%%%%%%%%%%%%%%%%%%%%%%%%%%%%%%%%%%%%%%%%%%%%%%%%%%%%%%%%%%%%%%

The Morse theory of the Yang-Mills functional on the space of gauge equivalence classes of
unitary connections on a hermitian vector bundle over a Riemann surface was introduced in
the seminal paper of Atiyah and Bott  \cite{AB}.  Further properties of the gradient flow of
the functional were obtained in \cite{D} and  \cite{R}. The minimal critical set can be
identified on the one hand with conjugacy classes of (projective) unitary representations
of the fundamental group of the surface via the holonomy map, and on the other hand with
the moduli space of semistable holomorphic vector bundles via the theorem of Narasimhan and
Seshadri (cf.\   \cite{NS, Do1}), and the analysis of \cite{AB} shows that the Yang-Mills flow
can be used effectively to study the topology of this space.

In \cite{DW2}  we studied the behavior of the moduli space of vector bundles as the
conformal structure on the Riemann surface degenerates.  A natural question to pose is
whether the Yang-Mills flow itself behaves in a reasonable fashion under similar
degenerations. In the case where the metric degeneration is to a \emph{cone metric}, one may
regard such a description as a non-linear version of the convergence of eigenvalues and
eigenfunctions of the Laplace operator obtained in \cite{JW1}.

The specific problem we consider is the following: let $X$ denote a compact surface
of genus $g\geq 2$, with a prescribed point $p$, and set $X^\ast=X\setminus\{p\}$.  We
define
\begin{equation}  \label{E:Ra}
\Ra=\Hom_\alpha\left(\pi_1(X^\ast), SU(2)\right)\bigr/ SU(2)
\end{equation} where the subscript $\alpha$ indicates that the holonomy of the
representation around $p$ is conjugate to the diagonal matrix with entries $e^{\pm 2\pi
i\alpha}$.   We take $0\leq\alpha\leq1/2$.  Thus, $\R_0$ is identified with conjugacy classes of
$SU(2)$ representations of the fundamental group of the closed surface $X$.
Fixing a conformal structure $[\sigma^\ast]$ on $(X,\{p\})$ we can define a smooth,
surjective map (the ``Hecke Correspondence"; see \cite{DDW} and below)
$
\pisigma_\alpha : \Ra\lra \R_0
$ 
for $\alpha\neq 0, 1/2$, which is an $S^2$-fibration over the irreducibles. 
Roughly speaking, the map is defined as follows:  for a flat $SU(2)$ connection $A$ with
holonomy $\alpha$ about $p$, we act on $A$  by a singular \emph{complex} gauge transformation $g$ to
bring the connection $A_0=g(A)$ into the standard, trivial form $d$ on a neighborhood of $p$.  Thus,
$A_0$ may be regarded as a connection on the closed surface $X$ which is, however, no longer flat.
The map is then defined by using the Yang-Mills flow to obtain from $A_0$ a flat connection
$\pisigma_\alpha([A])$.
This definition clearly involves the choice of conformal structure in an important way, and
understanding this behavior  is the motivation for this paper.

The map $\pisigma_\alpha$ can be generalized to a  map $\pisigma_\ab : \Ra\to \Rb$
for any $0\leq \beta\leq \alpha<1/2$, and this will be a homeomorphism for $\beta\neq 0$
(see Sec.\  \ref{S:foliation}).  We obtain from this a 1-dimensional foliation
$\foliation$ of
$\R=\bigcup_{0<\alpha<1/2}\Ra$.  A similar question to the one posed above is the dependence of
this foliation on $[\sigma^\ast]$.  For example, while it may be intuitively clear that
$\foliation$ varies continuously with $[\sigma^\ast]$, a differentiable structure is less
obvious.  Our first result is thus an explicit determination of variational formulas governing
$\foliation$.  As a consequence, we prove that the foliation is $C^1$.  The second result is a
description of the behavior when the conformal structure degenerates.  Here we show that for a
dense set of boundary points on Teichm\"uller space, the foliations actually converge away from
some singularities.  

More precisely,
let $\T_{aug.}(g,1)$ denote the augmented Teichm\"uller
space (cf.\ \cite{A}) obtained by adding \emph{nodal} Riemann surfaces. 
These are obtained by collapsing a collection $\Phi$ of disjoint simple closed
 boundary incompressible curves on $X^\ast$.  In the topology of $\T_{aug.}(g,1)$,
nodal surfaces may be approached by  the ``pinching" degeneration familiar from the
Deligne-Mumford compactification of the moduli of curves.  Thus, there is a family of conformal
structures
$\sigma^\ast(\ell)$ on $X^\ast$ degenerating as $\ell\to 0$ to a  conformal structure
$\sigma^\ast(0)$ on $X^\ast\setminus{\Phi}$, where a tubular neighborhood of each curve in
$c\in\Phi$ is conformal to the coordinate axes in $\CBbb^2$, minus the origin.  The collapsed
curve $c$, corresponding to the origin, is called the node (see Sec.
\ref{S:conicdegeneration} for more details). 
To see the effect of pinching on  the representation variety, we first make the following:

\begin{Def}  \label{D:accidental} 
Given a system $\Phi$ of disjoint simple closed
boundary incompressible curves on $X^\ast$ and an irreducible representation $\rho:
\pi_1(X^\ast)\to SU(2)$, we say that $\rho$ is \emph{accidentally reducible} with respect to
$\Phi$ if the restriction of
$\rho$ to the fundamental group of any component of $X^\ast\setminus\Phi$ is
reducible. We denote by $\R^\Phi\subset\R$ the closed subspace of conjugacy classes of
accidentally reducible representations.
\end{Def}

\noindent In Sec.\ \ref{S:proof} we prove the:

\begin{MainThm}   For each equivalence class $\sigma^\ast(0)$ of nodal conformal structures  on
$X^\ast\setminus{\Phi}$ there exists a
smooth 1-dimensional foliation ${\mathcal F}^{[\sigma^\ast(0)]} \subset
\R\setminus \R^\Phi$ such that for all paths $[\sigma^\ast(\ell)]\to[\sigma^\ast(0)]$ in
$\T_{aug.}(g,1)$ as above, 
$
\foliationtau\lra{\mathcal F}^{[\sigma^\ast(0)]}
$ on $\R\setminus \R^\Phi$ in the Hausdorff sense.
\end{MainThm}

This result is in contrast to the algebraic situation:
for each conformal structure $[\sigma^\ast]$ one can identify the space $\Ra$ with
the moduli space of parabolic stable bundles $\Ma^{[\sigma^\ast]}$ (cf.\ \cite{MS}).  This 
universalizes over $\T(g,1)$
to define a holomorphic fibration $\widetilde\Ma$.  Furthermore, the Mehta-Seshadri
Theorem defines an identification of $\widetilde\Ma$ with the trivial fibration $\widetilde
\Ra=\T(g,1)\times\Ra$.  In this setting, the map 
$
\pisigma_\alpha:\Ma^{[\sigma^\ast]}\lra\M
$
is the elementary transformation at $p$; for $0<\beta<\alpha<1/2$, $\pisigma_\ab$ is simply the
identity.
 Using algebro-geometric methods, it is possible to compactify
$\widetilde\Ma$ and $\widetilde\M$ over the Deligne-Mumford compactification of $\M(g,1)$ by
adding the appropriate moduli space of torsion-free sheaves on nodal curves.  Furthermore, one can
show that the maps $\pisigma_\alpha$ extend holomorphically over the compactification.  This is
similar in spirit to the degeneration used in \cite{DW2}.
But the map defined in the Main Theorem differs significantly from this
algebraic compactification in the sense that there does not appear to be a version of the
Narasimhan-Seshadri theorem relating the moduli space of torsion-free sheaves on a nodal curve to
the subspace of $\Ra$ mentioned above.  

We make two further comments on the proof:
(i) In the case where the Yang-Mills flow starting from a connection $[A]$ converges to an
irreducible flat connection $[B]$, one can find a complex gauge transformation $g$ such that
$[B]=g[A]$.
  This is not true, however, if $[B]$ is reducible.  As $[B]$ moves closer
to the reducibles one expects to lose control of the $C^0$ bound on $g$.  Reducibility of
$[B]$ can be detected via the existence of a kernel for the associated Laplace operator
$\Delta_B$ acting on the traceless, skew-hermitian endomorphisms of $E$.  As $[B]$
approaches the reducibles, the first eigenvalue of this operator goes to zero, and one
might therefore look for an explicit
$C^0$ bound on $g$ which depends on $[B]$ only through the first eigenvalue of $\Delta_B$. 
Such an estimate exists and will be  used in Sec.\  \ref{S:proof}.  As the surface degenerates to a
cone metric we have convergence of eigenvalues.  The
$C^0$ bound on $g$ mentioned in the previous paragraph will persist, provided that the
limiting eigenvalues are non-zero.  This is the motivation for the definition of accidental
reducibles (Def.\  \ref{D:accidental}).

(ii)  By a theorem of Wolpert \cite{Wo} the Strebel boundary is at finite Weil-Petersson
distance from the base point $[\sigma_0]$ in Teichm\"uller space.  From the structure of the
first variation formula and the bounds mentioned above, one expects
$\pi_\ab^{[\sigma^\ast(\ell)]}[A]$ to have finite length, and hence to converge, for  paths
$[\sigma^\ast(\ell)]\to[\sigma^\ast(0)]$.  For other points in the boundary we have no such
conclusion, as they are at infinite Weil-Petersson distance.  Though our proof will be along more
direct lines, this is the basic explanation for the convergence in the Main Theorem, and for why an
extension of $\pi_\ab^{[\sigma^\ast]}$  to the full Thurston boundary remains a difficult
question.  An interesting technical point in the proof below, however, is a restriction on the cone
angle of the degenerate metric required to prove the convergence.

This paper is organized as follows:  in Sec.\  \ref{S:gaugetheory} we review Donaldson's
approach to the Yang-Mills flow and R\aa de's version of L. Simon's estimate.  We discuss
Simpson's flow for singular metrics, and we show that the flow at infinity of a singular
connection preserves the conjugacy class of the holonomy about the singularity  (see
Cor.\  \ref{C:holonomy}).  We also interpret this construction by passing to branched
covers.  All this permits a definition of the foliation $\foliation$.  Sec.\  \ref{S:diff}
contains the proof of
the first variational formula Thm.\  \ref{T:firstvariation} for the action of the complex gauge
group. In particular, we show that the first variation for a path $[g_\varepsilon A_\varepsilon]$
of gauge equivalence classes of flat connections is independent of the derivative $\dot g$.  This
may be regarded as a kind of analogue of Ahlfors' result for the first variation of the hyperbolic
area element for quasi-conformal maps. In Sec.\ 
\ref{S:eigenvalue} we prove estimates for the eigenvalue problems for sections of vector
bundles and degenerations to cone metrics.    Finally, the Main Theorem is proved in Sec.\
\ref{S:proof}.

A word concerning notation:  if $\sigma$ is a Riemannian metric on a surface $X$,  then
integrals over $X$, unless otherwise specified, will be assumed to be taken with respect to
$\sigma$.  If $\sigma$ is a conformal metric and $z$ a conformal coordinate, we sometimes
write $\sigma(z)|dz|^2$ for the area form.   If $f$ is  a function on $X$, or more
generally a section of a vector bundle $V$ with a fiber metric $H$, and
$\alpha\geq 1$, we set 
$$
\Vert f\Vert_\alpha =\left\{ \int_X |f|^\alpha\right\}^{1/\alpha}
$$ In the above, if $f$ is a section, then
 $|f|=|f|_H$ involves the fiber metric $H$. When we want to emphasize the choice of
metrics, we will write $\Vert f\Vert_{\alpha;\sigma}$, or $\Vert f\Vert_{\alpha;\sigma,
H}$, or perhaps even $\Vert f\Vert_{\alpha; H}$ when there is no risk of confusion of the
two metrics. Other Sobolev norms will be denoted with an explicit subscript, e.g.
$\Vert f\Vert_{L^2_1(\sigma)}$, with $\Vert f\Vert_\infty$ (resp. $\Vert
f\Vert_{\infty, H}$)  for the $L^\infty$ norm of $|f|$ (resp. $|f|_H$).  Finally, we also use the
following abbreviation for $2\times 2$ diagonal matrices:
$$
\diag(\lambda_1,\lambda_2)=\left(\begin{matrix} \lambda_1 &0\\ 0&\lambda_2\end{matrix}\right)
\ .
$$

\medskip
\noindent \emph{Acknowledgements.}  We would like to thank Peter Li and Karen Uhlenbeck for
discussions, and the Max-Planck Institute in Bonn and the MSRI for their generous
hospitality.  Thanks also to G. Berger for suggesting Lemma \ref{L:grouptheory}.

%%%%%%%%%%%%%%%%%%%%%%%%%%%%%%%%%%%%%%%%%%%%%%%%%%%%%%%%%%%%%%%%%%%%%%%%%%%%%

\section{Gauge Theory}       \label{S:gaugetheory}

%%%%%%%%%%%%%%%%%%%%%%%%%%%%%%%%%%%%%%%%%%%%%%%%%%%%%%%%%%%%%%%%%%%%%%%%%%%%%

%%%%%%%%%%%%%%%%%%%%%%%%%%%%%%%%%%%%%%%%%%%%%%%%%%%%%%%%%%

\subsection{Review of the Yang-Mills Flow}   \label{S:ymflow}

%%%%%%%%%%%%%%%%%%%%%%%%%%%%%%%%%%%%%%%%%%%%%%%%%%%%%%%%%%

Let $E$ be a rank $2$ hermitian vector bundle on a closed, compact  surface $X$.  The
Yang-Mills flow is given by the equation:
\begin{equation} \label{E:ymflow}
\frac{dA(t)}{dt} +  D_{A(t)}^\ast F_{A(t)} =0\; ,
\end{equation} 
where $A(t)$ is a time dependent connection, $F_{A}$ is the curvature, and
$D^\ast_A=-\ast D_{A}\ast$ is the $L^2$ formal adjoint of   covariant differentiation
$D_A$ with respect to the connection
$A$.  Eq.\ (\ref{E:ymflow})  depends on a choice of a Riemannian metric on $X$ through the
Hodge
$\ast$, and this dependence will be the main theme of the paper.

Let $\A$ denote the space of unitary connections on $E$ and $\G$ the group of unitary gauge
transformations.  Donaldson proved in \cite{Do2} that the Yang-Mills flow equations
(\ref{E:ymflow})  with initial condition:
\begin{equation} \label{E:initialcondition} A(0)=A_0\ ,
\end{equation} 
have a unique solution for all time in $A/\G$.   Subsequently, R\aa de was
able to prove that the initial value problem (\ref{E:ymflow})-(\ref{E:initialcondition})
has a unique solution for all time in $\A$ \emph{before} we mod out by the gauge group.  We will
not need R\aa de's existence result in this paper, although some of his estimates on the
asymptotics of the flow as $t\to
\infty$ will be important for our arguments. 

Eq.\ (\ref{E:ymflow}) is the $L^2$-gradient flow for the Yang-Mills functional $YM(A)=\Vert
F_A\Vert_2^2$, and for a solution $A(t)$ of (\ref{E:ymflow}),
\begin{equation}  \label{E:ymgradient}
\frac{d}{dt}YM(A(t)) = -\nabla YM(A(t)) = -\Vert D^\ast_{A(t)} F_{A(t)}\Vert^2_2\ .
\end{equation} 
The critical points of $YM(A)$ are solutions of the Yang-Mills equations
$D^\ast_A F_A=0$, and the minima are given by the Hermitian-Yang-Mills (or projectively
flat) equation $\ast F_A=\mu\I$ for a constant $\mu$.  

Donaldson's approach to solving eq.'s (\ref{E:ymflow})-(\ref{E:initialcondition})
up to gauge is to solve instead the non-linear heat equation:
\begin{equation}  \label{E:metricflow} H^{-1}(t)\frac{dH(t)}{dt}=-\sqrt{-1}\ast F_{\dbar_E,
H(t)}\ ,
\end{equation} for a family of metrics $H(t)$ with initial conditions:
\begin{equation}  \label{E:metricinitialcondition} H(0)=H_0\ .
\end{equation} 
In the above, $F_{\dbar_E, H(t)}$ denotes the curvature of the unique
connection compatible with holomorphic structure $\dbar_E$ on $E$ and unitary with respect
to the metric $H(t)$.  

The two systems yield the same solution in $\A/\G$, for if $A(t)=g(t)A_0$ is a solution
of (\ref{E:ymflow}), then $H(t)=g(t)g(t)^\ast H_0$ is a solution to (\ref{E:metricflow});
conversely, if $H(t)=h(t)H_0$ is a solution to (\ref{E:metricflow}), then
$A(t)=h^{1/2}(t)A_0$ is real gauge equivalent to a solution of (\ref{E:ymflow}).  See
\cite{Do2} for more details.
Also notice that it is easy to factor out the trace part of the connection and gauge
transformations in eq.'s (\ref{E:ymflow})-(\ref{E:initialcondition}) and
(\ref{E:metricflow})-(\ref{E:metricinitialcondition}).  Therefore, in the following we
shall assume that the solutions to the above equations all preserve determinants.

We now review the analogue of L. Simon's result, due to R\aa de, concerning the
asymptotics of the solutions to (\ref{E:ymflow})-(\ref{E:initialcondition}):
\begin{Prop}[{\cite[Prop.\  7.2]{R}}] \label{P:rade} Let $A_\infty$ be an
\emph{irreducible} Yang-Mills connection.  There exist constants $\varepsilon_1$, $c>0$ such
that for any $A$ satisfying $\Vert A-A_\infty\Vert_{L^4}\leq\varepsilon_1$, we have:
\begin{equation}  \label{E:rade}
\Vert D_A^\ast F_A\Vert_2\geq c\left| YM(A)-YM(A_\infty)\right|^{1/2}\ .
\end{equation}
\end{Prop}

\noindent Let $\ad E$ denote the bundle of skew-hermitian endomorphisms of $E$, and let
$(\ad E)_0$ denote the subbundle of traceless ones.  Then the constants $\varepsilon_1$, $c$
 depend only on the first eigenvalue of the Laplacian associated to
$A_\infty$ acting on  sections of $(\ad E)_0$ and on the  constant governing the inclusion
$L^2_1\hookrightarrow L^4$ (notice that by Kato's inequality this is essentially the Sobolev
constant ${\mathfrak s}_1$ for functions -- see (\ref{E:df}) below).

The initial value problem (\ref{E:metricflow})-(\ref{E:metricinitialcondition}) has
solutions over non-compact  surfaces as well. More precisely, let $X$ be a compact Riemann
surface
 as before, choose ``punctures" $p_1,\ldots, p_k \in X$,  and let $X^\prime=X\setminus\{
p_1,\ldots,p_k\}$.  Fix a metric on $X^\prime$ whose expression in terms of a conformal
coordinate $z$ on $X$ centered at any one of the $p_i$'s has the form
$ds^2=\sigma(z)|dz|^2$, with:
\begin{equation}  \label{E:simpsonmetric}
\int_X |\sigma(z)|^p |dz|^2 < \infty\ ,
\end{equation} for some $p>1$ (cf.\   \cite[Prop.\  2.4]{Si1}).  Let $E$ be a holomorphic
bundle on $X$.
 Given a hermitian metric $H_0$ on $E$ with $\Vert \ast
F_{\dbar_E,H_0}\Vert_{\infty} < \infty$, we define:
$$
\deg(E,H_0)=\sqrt{-1}\int_X \tr\left(\ast F_{\dbar_E, H_0}\right)\ .
$$ 
Then $E$ is said to be $H_0$-\emph{stable} (resp. \emph{semistable}) if for any proper
holomorphic subbundle
$F$ of
$E$ we have:
$$
\frac{\deg(F,H_0)}{\rk F} < \frac{\deg(E,H_0)}{\rk E}\qquad(\text{resp.}\ \leq)\ .
$$ 
 Simpson proved the following:

\begin{Thm}[{cf.\   \cite{Si1}, Prop.\  6.6 and the proof of Thm.\  1}] 
\label{T:simpson} 
{\rm (i) } Given a holomorphic vector bundle $E$ on $X^\prime$ with
hermitian metric $H_0$ such that 
$\Vert \ast F_{\dbar_E, H_0}\Vert_{\infty, H_0} < \infty$, there exists a unique solution
$H(t)=h(t)H_0$ to (\ref{E:metricflow})-(\ref{E:metricinitialcondition}) with constant
determinant and having  the
 property that  for any finite $T>0$,
\begin{equation} \label{E:hbound}
\sup_{T\geq t\geq 0}\Vert h(t)\Vert_{\infty;H_0} < \infty\ .
\end{equation}
\par\noindent  {\rm (ii)} If, in addition, $E$ is $H_0$-stable, then (\ref{E:hbound}) holds
uniformly for all $T$.  Furthermore, $H(t)=h(t)H_0$ converges weakly in $L^p_{2, loc.}$ to a
solution
$H(\infty)=h(\infty)H_0$ of the Hermitian-Yang-Mills equation with $\Vert
h(\infty)\Vert_{\infty;H_0} < \infty$.
\end{Thm}

This result has important consequences which we will need in the proof of our main
theorem:  given a holomorphic bundle $E\to X^\ast$ as above and a puncture $p\in X\setminus
X^\ast$, choose a local holomorphic coordinate $z: U\to\Delta$ centered at $p$ and defined
on a neighborhood $U$, and a local holomorphic frame $\{f_1, f_2\}$ over $U$ such that
$H_0$ is in the standard form $\diag(|z|^{2\alpha}, |z|^{-2\alpha})$.  The identification is
chosen so that a unitary frame $\{e_1,e_2\}$ is given by $e_1=|z|^{-\alpha}f_1$,
$e_2=|z|^{\alpha}f_2$, and the hermitian connection $D_0$ with respect to this frame is in
the form
$D_0= d+\diag(i\alpha, -i\alpha)d\theta$.

\begin{Prop}  \label{P:diagonal}
The gauge
transformation $h(\infty)$ from Thm.\  \ref{T:simpson}, Part (ii), is independent of
the choice of conformal metric satisfying (\ref{E:simpsonmetric}).  Furthermore, if $\alpha< 1/2$
then $h(\infty)$ extends continuously at
$p$, and 
$h(\infty)(p)$ is diagonal with respect to the frame
$\{e_1,e_2\}$.
\end{Prop}

\begin{proof} For the first statement, note that the action of the complex gauge group
$\Gc$ on $\A$ is independent of the conformal factor.  Therefore, the $C^0$ bound from
Thm.\  \ref{T:simpson}, (ii), and the argument in \cite{Do1} prove uniqueness.
To show that $h(\infty)$ extends if $\alpha< 1/2$,   denote
by
$D_0$ the hermitian connection on
$E$ associated to
$H_0$.  Let
$g_1$ be a singular gauge transformation of the form $g_1=\diag(|z|^{-\alpha},
|z|^{\alpha})$ near $p$.  Then
$g_1(D_0)=d$ near $p$.  Furthermore, there exists $0\leq\beta\leq 1/2$ and a real gauge
transformation $\ell$ with $\det\ell=1$  such that if $g_2=\diag(|z|^{-\beta}, |z|^{\beta})$ near
$p$, then
$g_2\ell h^{1/2}(\infty)(D_0)=d$ (cf.\   \cite[Lemma 2.7]{DW1}).  It follows that $g_2\ell
h^{1/2}(\infty)g_1^{-1}(d)=d$, hence $g_2\ell h^{1/2}(\infty)g_1^{-1}$ is holomorphic on the
punctured disk $\Delta^\ast$.  On the other hand, since $|\alpha+\beta|<1$, this matrix
cannot have a pole at $p$, and it therefore extends continuously. Thus, we may write:
$$ g_2\ell h^{1/2}(\infty)g_1^{-1}=\left(\begin{matrix} a & b \\ c & d \end{matrix}\right)\ ,
$$ where $a,b,c,d$ are holomorphic in $\Delta$.  Then:
$$
\ell h^{1/2}(\infty)= g_2^{-1}\left(\begin{matrix} a & b \\ c & d \end{matrix}\right)g_1=
\left(\begin{matrix} |z|^{\beta-\alpha}a & |z|^{\alpha+\beta}b \\
|z|^{-(\alpha+\beta)}c & |z|^{\alpha-\beta}d \end{matrix}\right)\ .
$$ 
Suppose that $\alpha <\beta$.  Since the entries of the last matrix are bounded,
we must have $c(0)=0$ and $d(0)=0$, which contradicts $\det h^{1/2}(\infty)=1$  (as mentioned above
we always fix the determinant by projecting away the trace part of the connection; see \cite{Si1}). 
A similar argument holds for $\alpha>\beta$.  Hence, $\alpha=\beta$.  Finally, since $b$ and $c$ are
holomorphic, $c(0)=0$, and $\ell h^{1/2}(\infty)$ is diagonal at $p$. In particular, $(\ell
h^{1/2}(\infty))^\ast(\ell h^{1/2}(\infty))(p)=h(\infty)(p) $ is diagonal at $p$.
\end{proof}

\begin{Cor}   \label{C:holonomy} If $E$ is a holomorphic bundle on $X^\ast$ which is
$H_0$-stable, then the Yang-Mills flow at infinity preserves the conjugacy class of the
holonomy around the punctures.
\end{Cor}

\begin{proof}
The case $0<\alpha<1/2$ follows from the proof of Prop.\  \ref{P:diagonal} above.  For
$\alpha=0$, the metric $H_0$ is smooth at $p$ and the flow
(\ref{E:metricflow})-(\ref{E:metricinitialcondition}) extends smoothly on $X$.  It follows that
the holonomy of the limit remains trivial.  The case $\alpha=1/2$ also follows from the proof of 
Prop.\  \ref{P:diagonal}, for if $\beta<1/2$ then we again have $|\alpha+\beta|<1$, and the same
argument gives a contradiction.
\end{proof}

\begin{Rem}  \label{R:stability} 
It follows by the arguments in \cite{Si2} that if $H_0$ is
of the form $\diag(|z|^{2\alpha}, |z|^{-2\alpha})$ with respect to a holomorphic frame
$\{f_1,f_2\}$ near $p$, then $H_0$-stability coincides with Seshadri's parabolic stability
with respect to the weights $\{\alpha, -\alpha\}$.  For more details we refer to
\cite{DW1,DW2}.
\end{Rem}

%%%%%%%%%%%%%%%%%%%%%%%%%%%%%%%%%%%%%%%%%%%%%%%%%%%%%%%%%%

\subsection{Representation Varieties and Branched Covers}  \label{S:covers}

%%%%%%%%%%%%%%%%%%%%%%%%%%%%%%%%%%%%%%%%%%%%%%%%%%%%%%%%%%

Let $X$ be a compact surface, $p\in X$, and $X^\ast=X\setminus\{p\}$.
We denote by $\R(X^\ast)$ the space of conjugacy classes of $SU(2)$ representations of the
free group $\pi_1(X^\ast)$.  We may also identify $\R(X^\ast)$ with the space of gauge
equivalence classes of flat connections on a trivial rank two hermitian bundle $E$ on
$X^\ast$. Given a real number $0\leq \alpha\leq 1/2$, we denote by $\Ra$ the subspace of
flat connections on $E$ with holonomy matrix conjugate to $\exp\left\{ 2\pi
i\diag(\alpha,-\alpha)\right\}$ around the puncture $p$.  Then $\R_0$ is naturally identified with
the space
$\R(X)$ of equivalence classes of flat connections on the trivial bundle over the compact
surface $X$.  Notice that for $\alpha\neq 0$, $\Ra$ consists entirely of irreducible
representations.  For future reference, we let $\R(X)_{irr}\subset\R(X)$ denote the open
set of irreducible representations of $\pi_1(X)$.

In this section, we will sketch how to reduce certain analytical questions for the
Yang-Mills flow on $E$ over $X^\ast$ to the flow  on a bundle over a branched cover of the
compact surface
$X$. 
Let $\beta=k/n$ where $k, n$ are positive coprime integers.  Consider a regular, cyclic,
$n$-fold, holomorphic branched cover $\widehat X$ of $X$ with $p$ in the ramification divisor $B$. 
Let
$ q:\widehat X\to X$ be the covering map, $\hat p=q^{-1}(p)$, $\widehat B=q^{-1}(B)$, 
$\widehat X^\ast=\widehat X\setminus\{\hat p\}$, $\widehat U=q^{-1}(U)$, and $\widehat
U^\ast=\widehat U\cap \widehat X^\ast$.  Let $E=X^\ast\times\CBbb^2$ be the trivial rank 2
vector bundle on $X^\ast$.  We construct a bundle $\widehat E$ over $\widehat X$ by gluing
$q^\ast(E)$ on $\widehat X\setminus\widehat U$ with $\widehat U\times\CBbb^2$ via the gauge
transformation:
$$
\hat s : \widehat U^\ast\lra SU(2) \quad :\quad
\hat s(w)= \hat s(\hat r,\hat\theta)=
\left(\begin{matrix} e^{-ik\hat\theta} & 0\\ 0 & e^{-ik\hat\theta} \end{matrix}\right)\ .
$$
Since $\det \hat s=1$, it follows that $\deg\widehat E=0$; hence, $\widehat E$ is isomorphic
to the trivial bundle. Although there is no natural global trivialization, there is a
 trivialization  of $\widehat E\bigr|_{\widehat X^\ast}=q^\ast E\bigr|_{\widehat
X^\ast}$ induced by the one on $E$.

On $E$ we fix the trivial hermitian metric $H_0$.  Also, since $\hat s$ is unitary, the
trivial metrics on $q^\ast E$ and $\widehat U\times\CBbb^2$ glue together to define a
hermitian metric $\widehat H_0$ on $\widehat E$.  Let $\Ab$ denote the space of unitary
connections of $E$ on $X^\ast$, flat in a neighborhood of $p$ with holonomy conjugate to
$\exp\left\{ 2\pi i\diag(\beta,-\beta)\right\}$.  Given $A\in \Ab$, choose a real gauge
transformation
$g$ such that:
$$
g(A)\bigr|_{U}= d+\left(\begin{matrix} i\beta & 0\\ 0 & -i\beta
\end{matrix}\right)d\theta\ .
$$
It follows that:
$$q^\ast g(A)\bigr|_{\widehat U^\ast}=d+\left(\begin{matrix} ik & 0\\ 0 & -ik
\end{matrix}\right)d\hat\theta\ .
$$
By gluing $q^\ast g(A)$ with the trivial connection $d_{\widehat U}$ via $\hat s$, we
obtain a unitary connection $\widehat A$ on $\widehat E$.  Let $\widehat\A$ denote the
space of unitary connections on $\widehat E$.  Define:
$$
\hat q:\Ab\lra\widehat\A\quad ,\quad \hat  q(A)=q^\ast(g^{-1})\widehat A\ ,
$$
where in the above, $q^\ast(g)=g\circ q$.  It is easily checked that $\hat q$ is
well-defined and real gauge equivariant.  In particular, it induces a map
$
\hat q:\Rb\to\widehat\R=\R(\widehat X)=\widehat\A_{flat}/\widehat\G
$,
where $\widehat\A_{flat}$ are the flat connections on $\widehat E$ and $\widehat\G$ is the
real gauge group.  We note the following:

\begin{Prop}   \label{P:branchedreducibles}
If $\beta=k/n$ with $n$ odd, then $\hat q(\Rb)\subset \widehat R_{irr}$.  Furthermore, given a
collection $\Phi$ of simple closed boundary incompressible curves in $X^\ast$, mutually
disjoint and disjoint from $B$, $\widehat {\Phi}=q^{-1}({\Phi})$, and $[A]\in\Rb$
which is not accidentally reducible with respect to $\Phi$ (see Def.\ \ref{D:accidental}),
then
$\hat q([A])$ is not accidentally reducible with respect to $\widehat {\Phi}$.
\end{Prop}

\begin{proof}  We first prove:

\begin{Lem}  \label{L:grouptheory}
Consider an exact sequence of groups $1\to H\to G\to Q\to 1$, and suppose that $Q$ is
abelian with no index 2 subgroup.  Then the restriction to $H$ of any irreducible $SU(2)$
representation of $G$ is also irreducible.
\end{Lem}

\begin{proof}
Suppose, to the contrary, that there is an irreducible  $\rho:G\to SU(2)$
which is reducible on $H$.  Then there is a maximal torus $\bf T$ in $SU(2)$ such that the
image of $H$ lies in $\bf T$.  Since $Q$ is abelian and $\rho$ is irreducible, the image of
$H$ cannot be contained entirely in the center of $SU(2)$.  This implies that the image of
$G$ lies in
$N({\bf T})$, the normalizer of $\bf T$.  Now we have the exact sequence for the Weyl group
$1\to {\bf T}\to N({\bf T})\to \ZBbb/2\ZBbb\to 1$.  Since the image of $G$
is not contained in $\bf T$, we must have a surjection $Q\to\ZBbb/2\ZBbb\to 1$;
contradiction.
\end{proof}

\noindent  Continuing with the proof of Prop.\  \ref{P:branchedreducibles}: as remarked
above, a point $[A]\in\Rb$ gives a conjugacy class of irreducible representations of
$\pi_1(X^\ast)$ which, in turn, induce irreducible representations of $\pi_1(X\setminus B)$.  Now
$\hat q[A]$ is irreducible as a representation of $\pi_1(\widehat X)$ if and only if it induces
an irreducible representation of $\pi_1(\widehat X\setminus{\widehat B})$.  The first statement
then follows from Lemma \ref{L:grouptheory} by setting $H=\pi_1(\widehat X\setminus{\widehat B})$,
$G=\pi_1( X\setminus B)$, and $Q=\ZBbb/n\ZBbb$.  For the second statement, consider a connected
component
$Y$ of $X\setminus{\Phi}$.  If $Y\cap B\neq \emptyset$, then $\widehat Y=q^{-1}(Y)$ is a
connected component of $\widehat X\setminus\widehat{\Phi}$, and the argument is as above. 
If $Y\cap B=\emptyset$, then for each connected component $\widehat Y_i$ of $q^{-1}(Y)$,
$\pi_1(\widehat Y_i)\simeq\pi_1(Y)$, so $\hat q[A]$ is irreducible there as well.
\end{proof}

Now consider the effect of eq.'s (\ref{E:metricflow})-(\ref{E:metricinitialcondition})
under $\hat q$.   First we choose a conformal metric $\sigma$ on $X\setminus B$ with the
following property:  for every $b\in B$, express the map $q$ locally as $z=w^n$ for 
coordinates
$z$ on
 centered at
$b$ and $w$ centered at $\hat b$.  Then we assume that $\sigma(z)=1/n^2|z|^{2(1-1/n)}$. 
Notice that such a metric satisfies condition (\ref{E:simpsonmetric}).  This is an example
of a cone metric (see (\ref{E:conemetric})).  Let $\hat\sigma$ be the pull-back metric on
$\widehat X\setminus {\widehat B}$.  Then the condition on $\sigma$ implies that $\hat
\sigma$ extends to a smooth conformal metric on $\widehat X$.

Next, fix a connection $A_0$ on $E$ over
$X^\ast$ and let $\widehat A_0=\hat q(A_0)$.  Let $H(t)=h(t)H_0$ be a solution of 
(\ref{E:metricflow})-(\ref{E:metricinitialcondition}) on $E$ over $X\setminus B$, where the
holomorphic structure on
$E$ is defined by $A_0^{0,1}$.  Then $\widehat H(t)=q^\ast H(t)=\hat h(t)\widehat H_0$ is a
solution of the same equations on $\widehat X\setminus{\widehat B}$ with respect to the
holomorphic structure $q^\ast A_0^{0,1}$.  Since $q$ is smooth and $h(t)$ satisfies the
estimate (\ref{E:hbound}), the same is true for $\hat h(t)$.  Hence, by the uniqueness
properties of the flow, we obtain
\begin{Lem}  \label{L:branchedflow}  
The restriction of the flow $\widehat H(t)$ to $\widehat X\setminus{\widehat B}$ coincides
with
$q^\ast H(t)$.
\end{Lem} 

\noindent  This allows us to reduce estimates of $h(t)$ over $X^\ast$ to estimates of $\hat
h(t)$ over
$\widehat X$. For example:

\begin{Lem}  \label{L:branchedestimate}
Let $\Omega$ be an open set with compact closure in $X\setminus B$.  Then there exists a
 constant $C=C(\Omega)$ such that:
$$
\Vert h(t)\Vert_{L^p_k(\Omega)}\leq  C\, \Vert \hat h(t)\Vert_{L^p_k(\widehat\Omega)}\ .
$$
In the above, $h(t)$ is the solution to
(\ref{E:metricflow})-(\ref{E:metricinitialcondition}) on $X^\ast$, and $\hat h(t)$ is the
solution over
$\widehat X$.
\end{Lem}

%%%%%%%%%%%%%%%%%%%%%%%%%%%%%%%%%%%%%%%%%%%%%%%%%%%%%%%%%%

\subsection{Definition of the Foliation}   \label{S:foliation}

%%%%%%%%%%%%%%%%%%%%%%%%%%%%%%%%%%%%%%%%%%%%%%%%%%%%%%%%%%

Let $X$, $X^\ast$ be as above.  
Given $0\leq \beta <\alpha < 1/2$, or $0<\alpha<\beta\leq 1/2$, and
$[\sigma^\ast]\in\T(g,1)$, we will now give a rigorous definition of the twist maps
$\pisigma_{\ab}:\Ra\to\Rb$. The conformal structure on $X^\ast$ extends to $X$.  Let
$\sigma$ denote a choice of smooth metric on $X$ compatible with this conformal structure. 
Let
$z:U\to\Delta$ be a holomorphic coordinate centered at $p$, and let
$\varphi_1$ be a smooth cut-off function supported in $\Delta_{1/3}$ and identically equal
to $1$ on
$\Delta_{1/6}$.  We define a singular complex gauge transformation of $E$ by setting:
\begin{equation} \label{E:g} g_{\ab}(\xi)=
\exp\left\{\varphi_1(z(\xi))
          \left(    {\begin{matrix} (\alpha-\beta)\log|z(\xi)| & 0 \\
                         0 & (\beta-\alpha)\log|z(\xi)| 
                       \end{matrix}  }
          \right)
         \right\}, 
\end{equation} for $\xi\in U$ and extending by the identity elsewhere. We now come to the
definition of $\pisigma_\ab$:

\begin{Def} \label{D:foliation1}
Given $[A]\in\Ra$, choose a representative $A$ in the standard form
$\diag(i\alpha , -i\alpha)d\theta$  when pulled back to $\Delta$ via $z$.  Then $\widetilde
A=g_\ab(A)$ will be of the form $\diag(i\beta , -i\beta)d\theta$ over $\Delta_{1/6}$. 
Furthermore, it is easy to check (cf.\   \cite{MS,DW1}) that the holomorphic structure  induced by
$\widetilde A$ is parabolic stable, and hence by Remark \ref{R:stability} we can flow the
connection
$\widetilde A$ with respect to the metric $\sigma$ at infinite time.  Thm.\ 
\ref{T:simpson} guarantees that
$\widetilde A$ flows to a flat connection $\widetilde A(\infty)$ which, by Cor.\ 
\ref{C:holonomy} has holonomy conjugate to
$\exp\{2\pi i\diag(\beta,-\beta)\}$.  We set:
\begin{equation}  \label{E:pidefinition}
\pisigma_\ab[A]=[\widetilde A(\infty)]\ ,
\end{equation} the class of $\widetilde A(\infty)$ in $\Rb$.  
\end{Def}

It is straightforward to show that the definition of $\pisigma_\ab$ is independent of the
choices made, i.e.\ the coordinate $z$, the cutoff function $\varphi_1$, the lift of the
conformal structure, and the lift of $[A]$.  In other words, the only dependence of
$\pisigma_\ab$ is through the class $[\sigma^\ast]\in T(g,1)$.
Set $\R=\bigcup_{0<\alpha<1/2} \Ra \subset \R(X^\ast)$.  The next two lemmas are left to
the reader:

\begin{Lem}  \label{L:composition} For $0\leq\gamma<\beta<\alpha<1/2$ or
$0<\alpha<\beta<\gamma\leq 1/2$ and $[\sigma^\ast]\in\T(g,1)$,
$\pisigma_{\alpha\gamma}=\pisigma_{\beta\gamma}\circ\pisigma_\ab$.
\end{Lem}

\begin{Lem}  \label{L:picontinuity} 
Given $[\sigma^\ast]\in\T(g,1)$, $\alpha\in(0,1/2)$, and $[A]\in\Ra$, the map
$[0,1/2]\to\R:
\beta\mapsto\pisigma_\ab[A]$ is continuous.  Furthermore, its restriction to $(0,1/2)$ is
smooth.
\end{Lem}

\noindent For convenience, for $\beta=\alpha$ we define $\pisigma_\ab$ to be the identity.

\begin{Def}  \label{D:foliation2} Fix $[\sigma^\ast]\in T(g,1)$, $[A]\in\Ra$,
$0<\alpha<1/2$, and let:
$$
\foliation_{[A]}=\bigcup_{0<\beta<1/2}\pisigma_\ab[A]\quad , \quad
\foliation =\bigcup_{[A]\in\R^\ast}\foliation_{[A]}\ .
$$
\end{Def}

\noindent It follows from the above discussion that $\foliation$ is a smooth
1-dimensional foliation of
$\R^\ast=\R\setminus\hol^{-1}\{0,1/2\}$.

We now turn to the definition of the limiting foliation $\foliationphi$ discussed in the
Introduction. 
Let $\Phi$ denote a collection of disjoint simple closed curves on $X^\ast$, and let
$\sigma^\ast(0)$ be a conformal structure on the pinched surface $X_0^\ast=X^\ast\setminus\Phi$.
Let us elaborate on this (see Sec.\ \ref{S:conicdegeneration} for more details).  By definition,
the conformal structure on the pinched surface $X^\ast_0$ arises from a conformal structure
$\sigma^\ast$ on $X^\ast$ as follows:  for each $c\in \Phi$ there is a tubular, cylindrical
neighborhood $C$ of $c$ which is conformally equivalent, with respect to $\sigma^\ast$, to
the intersection in $\CBbb^2$ of a
neighborhood of the origin with the annulus $zw=\varepsilon$ for some non-zero complex number
$\varepsilon$.  In these coordinates, $c$ is the set where $|z|=|w|=|\varepsilon|^{1/2}$.  The
conformal structure $\sigma^\ast(0)$ then replaces $C$ with the pinched annulus $zw=0$. 

  Let
$\R^\Phi\subset\R$ be the set of
$\Phi$-accidental reducibles.  Given
$[A]\in
\Ra\setminus \R^\Phi$, let $A$ be a lift of
$[A]$ to a connection which has the standard form $d+\diag(i\gamma,-i\gamma)d\theta$
in a coordinate neighborhood $z$ of each $c\in\Phi$ as described above. Of course, the holonomy
$\gamma$ depends on the component $c$. By assumption, the restriction of
$A$ to any component of
$X_0^\ast$ is irreducible.  Hence, for $\varepsilon([A]) >0$ sufficiently small and
$|\beta-\alpha| < \varepsilon([A])$, the holomorphic structure associated to the twists
$g_\ab(A)$ on $X^\ast_0$ is parabolic stable for the choice of weights $\beta$ and
$\{\gamma_c\}_{c\in \Phi}$.

\begin{Def}  \label{D:infiniteflow} 
Set $A(0,\infty)$ to be the flow of $g_\ab(A)$ (at
infinite time) with respect to a metric on $X_0^\ast$ compatible with the conformal
structure and satisfying the condition (\ref{E:simpsonmetric}).  As before, the real gauge
equivalence class $[A(0,\infty)]$ of
$A(0,\infty)$ is independent of all the choices made.
\end{Def}

\begin{Thm}   \label{T:lift}
 There is a well-defined lift of $[A(0,\infty)]$ to an element in
$\Rb$.
\end{Thm}

\begin{proof}Let $U_\pm$ be coordinate neighborhoods $z$ and $w$
corresponding to a curve
$c\in\Phi$ as before,
and suppose the holonomy is $\gamma$.  Thus we have:
$$ A\bigr|_{U_\pm} = d \pm \left(\begin{matrix} i\gamma & 0 \\ 0 &
-i\gamma\end{matrix}\right)d\theta_\pm\ ,
$$ with respect to unitary frames $e_1^\pm, e_2^\pm$ of $E|_{U_\pm}$.  By Cor.\ 
\ref{C:holonomy} there are real gauge transformations $g_\pm$ of $E|_{U_\pm}$ such that:

\begin{equation}  \label{E:connection} g_\pm\left(A(0,\infty)\right)\bigr|_{U_\pm} = d \pm
\left(\begin{matrix} i\gamma & 0 \\ 0 & -i\gamma\end{matrix}\right)d\theta_\pm\ .
\end{equation}

\noindent The $g_\pm$ may be extended to a global real gauge transformation $g$ of
$E|_{X_0^\ast}$ which is the identity away from a small neighborhood of $U_\pm$.  By using
the identification of the frame
$\{e_1^+,e_2^+\}$ with
$\{e_1^-,e_2^-\}$, (\ref{E:connection}) implies that the pull-back of $g(A(0,\infty))$ to
$X^\ast\setminus c$ extends smoothly over $c$.  By repeating the above for every curve
$c\in\Phi$ we obtain a flat connection on $X^\ast$ with the correct holonomy. 
The resulting conjugacy class is unique, because a conjugacy class of connections on
$E|_{X_0^\ast}$ together with gluing data determine a unique conjugacy class of connections
on
$E|_{X^\ast}$.
\end{proof}

The construction above may be summarized as follows:  the connection $A$ on $X^\ast$ determines a
flat connection on $X^\ast\setminus\Phi$, along with ``gluing parameters" across the curves
$c\in\Phi$.  The twisted connection $g_\ab(A)$ has a small amount of curvature in the component $Y$
of $X^\ast\setminus \Phi$ containing $\{p\}$, and it agrees with $A$ on all the other components. 
The flow $A(0,\infty)$ then runs the Simpson flow on $Y$, applied to $g_\ab(A)$, and leaves the
connection on the other components fixed.  Since the holonomies around the curves $c\in\Phi$
bounding $Y$ remain unchanged for the flow at infinite time (Cor.\ \ref{C:holonomy}), we may use the
original gluing parameters to reconstruct a flat connection on $X^\ast$.

\begin{Def}  \label{D:infinitefoliation} For each $[A]\in\Ra$ and $\beta$ satisfying
$|\alpha-\beta|<\varepsilon([A])$ as above, define
$\pi_\ab^{[\sigma^\ast(0)]}[A]$ to be the lift of $[A(0,\infty)]$ described in Thm.\ 
\ref{T:lift}.  Furthermore, we set:
$$
\foliationphi_{[A]}=\bigcup_{\alpha-\varepsilon([A])<\beta<\alpha+\varepsilon([A])}
\pi_{\ab}^{\Phi}[A]\ .
$$ It is straightforward to check a composition rule as in Lemma \ref{L:composition}.  We
may therefore define:
$$
\foliationphi = \bigcup_{[A]\in\R} \foliationphi_{[A]}\ .
$$ 
Then $\foliationphi$ is a smooth foliation of $\R\setminus\Phi^\Phi$.
\end{Def}

%%%%%%%%%%%%%%%%%%%%%%%%%%%%%%%%%%%%%%%%%%%%%%%%%%%%%%%%%%

\section{Differentiability of the Foliation}   \label{S:diff}

%%%%%%%%%%%%%%%%%%%%%%%%%%%%%%%%%%%%%%%%%%%%%%%%%%%%%%%%%%%

%%%%%%%%%%%%%%%%%%%%%%%%%%%%%%%%%%%%%%%%%%%%%%%%%%%%%%%%%%

\subsection{First Order Variational Formula}    \label{S:variationalformula}

%%%%%%%%%%%%%%%%%%%%%%%%%%%%%%%%%%%%%%%%%%%%%%%%%%%%%%%%%%

In describing the behavior of the Yang-Mills flow applied to a connection $A$ as the conformal
structure $\sigma$ on a closed surface varies, there are two considerations:  first, the complex
gauge transformation $g$ describing the flow, i.e.\ such that $g(A)$ is flat, will depend in a
complicated way on $\sigma$.  Second, while the real gauge group acts on the space of unitary
connections in a manner independent of the conformal structure, the complex gauge group does not.

Assume that we have fixed a
conformal structure on $X$.  Let
$\mu_{\varepsilon}$ be a differentiable family of Beltrami differentials on $X$ with $\mu_0=0$,
$\dot\mu_\varepsilon\bigr|_{\varepsilon=0}=\nu$.  Let
$g_{\varepsilon}$ be a differentiable family of complex gauge transformations on $E$ with $g_0=g$,
$\dot g_\varepsilon\bigr|_{\varepsilon=0}=\dot g$.  Finally, let $A_{\varepsilon}$ be
differentiable family of unitary connections on $E$ with
$A_0=A$, $\dot A_\varepsilon\bigr|_{\varepsilon=0}=\dot A$.  Set
$\gamma_{\varepsilon}=g_{\varepsilon} A_{\varepsilon}$, and $\dot \gamma=\dot
\gamma_\varepsilon\bigr|_{\varepsilon=0}$. We emphasize that the action of the complex
gauge group is with respect to the
$\mu_{\varepsilon}$-deformed conformal structure. Also,  we regard Beltrami differentials $\nu$ as
endomorphisms $\Omega^{1,0}(X)
\to\Omega^{0,1}(X)$, which extend to endomorphisms on forms with values in $E$.

\begin{Thm} \label{T:firstvariation}  
 Let $\mu_{\varepsilon}$, $g_{\varepsilon}$, $A_{\varepsilon}$, and
$\gamma_{\varepsilon}$ be as above.  Then:
$$
\dot\gamma^{0,1}=\dbar_{g(A)}\left(g^{-1}\dot g\right)-\nu\left[ \left(\partial_A
g^\ast\right)(g^\ast)^{-1}+(g^\ast)^{-1}\partial_A g^\ast\right]+g^{-1}\dot A^{0,1} g\ ,
$$
 where the $(0,1)$ part is taken with respect to the fixed conformal structure on $X$,
and $g^\ast$ denotes the fiberwise adjoint of $g$ with respect to the background metric $H_0$.
\end{Thm}

\begin{Rem}
If $\gamma_\varepsilon$ is a path of flat connections, then in $\A/\G$ a representative for the
tangent vector
$[\dot\gamma^{0,1}]$ may be taken to be harmonic.  Projecting to the harmonics, we see that the
first term in the expression above vanishes.  Hence, we conclude that the first variation of
the Yang-Mills flow is independent of $\dot g$.
\end{Rem}

\noindent \emph{Proof of Thm.\ \ref{T:firstvariation}}.
We first note how the Cauchy-Riemann operators deform:  fix a unitary connection $A$ on $E$ and
a Beltrami differential $\mu$.  
\begin{Lem}  \label{L:dbardeformation}
Let   
$\dbar_{A,\mu}:\Omega^0(E)\to\Omega^{0,1}(X_\mu, E)\subset\Omega^1(E)$
be the $\dbar$-operator on $X_\mu$ associated to $A$.  Then $\dbar_{A,\mu}$ is given by
$$
\dbar_{A,\mu}s=\frac{1}{1-|\mu|^2}\left(\dbar_A s-\mu\partial_A s +\bar\mu\dbar_A
s-|\mu|^2\partial_A s 
\right)\ ,
$$
for smooth sections $s\in\Omega^0(E)$.
\end{Lem}

\begin{proof}
By choosing local frames, the lemma follows from the corresponding statement for $\dbar$ acting
on functions.  Thus, let $f$ be a function, $z$ a local conformal coordinate on $X$, and
$w=w_{\varepsilon}$  solutions to the Beltrami equation $ w_{\bar z}=\mu^z_{\bar z}
w_z$. 
Here we have expressed $\mu=\mu_{\bar z}^z d\bar z\otimes(\partial/\partial z)$.  
 Write:
$$
df=f_z dz+ f_{\bar z} d\bar z=f_w dw+ f_{\bar w} d\bar w\ ,
$$
and use $dw=w_z(dz+\mu_{\bar z}^z d\bar z)$ to obtain:
$$
f_z=f_z w_z+f_{\bar w}\bar{w_z}\bar\mu_z^{\bar z}\ , \
f_{\bar z}=f_w w_z\mu_{\bar z}^z + f_{\bar w}\bar{w_z}\ .
$$
Multiplying the first equation by $\mu$ and subtracting the second, we have:
$$
f_{\bar w}\bar{w_z}(1-|\mu|^2)=f_{\bar z}-f_z\mu_{\bar z}^z\ .
$$
Now multiply through by $d\bar z+\bar\mu_{\bar z}^z dz$ to obtain:
$$
\dbar_\mu f=f_{\bar w}d\bar w=\frac{1}{1-|\mu|^2}\left(f_{\bar z}d\bar z-\mu_{\bar z}^z d\bar z f_z
+\bar\mu_{\bar z}^z dz f_{\bar z}-|\mu|^2 f_z dz
\right)\ .
$$
The result follows by observing that:
$$
\mu_{\bar z}^z d\bar z f_z=\mu_{\bar z}^z d\bar z\otimes(\partial/\partial
z)(f_z dz)=\mu\partial f\ , \
\bar\mu_{\bar z}^z d z f_{\bar z}=\bar\mu^{\bar z}_z dz\otimes(\partial/\partial
z)(f_{\bar z} d\bar z)=
\bar\mu\bar\partial f\ .
$$
\end{proof}
Continuing with the proof of Thm.\  \ref{T:firstvariation}: the action of the complex
gauge group is given by:
\begin{equation}\label{E:complexgaugegroupaction}
g(D_A)=d+ g^{-1}A^{0,1}g+g^\ast A^{1,0}(g^\ast)^{-1}+g^{-1}\dbar_\mu g-(\partial_\mu
g^\ast)(g^\ast)^{-1}\ ,
\end{equation}
where
\begin{equation} \label{E:Adeformation}
\begin{aligned}
A^{0,1}&= A_{\bar w}d\bar w=A_{\bar z}d\bar z-\varepsilon\nu(A_zdz)+\varepsilon\bar\nu(A_{\bar
z}d\bar z)+O(\varepsilon^2)\ ,\\
A^{1,0}&= A_{ w}d w=A_{ z}d z-\varepsilon\bar\nu(A_{\bar z}d\bar z)+\varepsilon\nu(A_{
z}dz)+O(\varepsilon^2)\ .
\end{aligned}
\end{equation}
From Lemma \ref{L:dbardeformation} we have:
\begin{align*}
\dbar_\mu g &=\dbar g+\varepsilon\left(\dbar\dot g-\nu(\partial g)+\bar\nu(\dbar
g)\right)+O(\varepsilon^2)\ ,\\
\partial_\mu g^\ast &=\partial g^\ast+\varepsilon\left(\partial\dot g^\ast-\bar\nu(\dbar
g^\ast)+\nu(\partial g^\ast)\right)+O(\varepsilon^2)\ .
\end{align*}
Applying this and (\ref{E:Adeformation}) to (\ref{E:complexgaugegroupaction}) yields the result.
\qed

%%%%%%%%%%%%%%%%%%%%%%%%%%%%%%%%%%%%%%%%%%%%%%%%%%%%%%%%%%

\subsection{Differentiability of the Twisted Connection}   \label{S:twistbound}

%%%%%%%%%%%%%%%%%%%%%%%%%%%%%%%%%%%%%%%%%%%%%%%%%%%%%%%%%%

The aim of this section is to prove the following:

\begin{Prop}  \label{P:critical}
 Let $\sigma^\ast_{\varepsilon}$ be  a continuously differentiable family of the metrics 
representing a path $[\sigma^\ast_{\varepsilon}]\subset\T(g,1)$. Then there is a family of singular
gauge transformations
$g_{\ab}^{\varepsilon}$ such that $g_{\ab}^{\varepsilon}(A)$ is a continuously differentiable path
of connections.
\end{Prop}

\begin{proof}
We differentiate at $\varepsilon=0$.   Let $\mu_{\varepsilon}$ be the differentiable
path of Beltrami differentials on
$X$ with $\mu_0=0$ and $\dot\mu_\varepsilon\bigr|_{\varepsilon=0}=\nu$, associated to
$\sigma^\ast_{\varepsilon}$.  Let $z:U\to\Delta$ be a local  coordinate on a neighborhood
$U$ of $p$, conformal with respect to the conformal structure determined by $\sigma_0^\ast$.
 Fix
$\varphi_0$ a smooth cut-off function supported in
$\Delta$ and identically 1 on $\Delta_{2/3}$.
We obtain  Beltrami differentials
$\tilde\mu$ on
$\CBbb$ by extending
$\varphi_0\mu$ by zero.  For $\varepsilon\geq 0$ small, we consider the solution
$w_{\varepsilon}$ to the Beltrami equation on $\CBbb$:
\begin{equation}  \label{E:beltrami2}
 w_{\bar z}=\tilde\mu^z_{\bar z}
w_z\ ,
\end{equation}
 normalized such that:
\begin{equation}  \label{E:normalization} w_{\varepsilon}(0)=0\ ,\quad w_{\varepsilon}(1)=1\
,\quad w_{\varepsilon}(\infty)=\infty\ .
\end{equation} 
Let $\dot w$ denote $\partial w_\varepsilon/\partial\varepsilon$ at $\varepsilon=0$.   We also set
$\tilde\nu=\dot{\tilde\mu}$.

To compute the derivative, we must take into account the change of frame.  The problem is local,
so suppose we have a fixed (trivial) hermitian rank 2
vector bundle over the complex plane with global unitary frame $e_\pm$ and singular
connection $A_\alpha$ such that:
\begin{equation}  \label{E:singularconnection} D_A e_\pm =\pm i\alpha d\theta\otimes e_\pm\ .
\end{equation} 
Let  $\theta_{\varepsilon}$ denote the theta coordinate of
$w_{\varepsilon}$.
\begin{Lem}  \label{L:frametwisting} Define a frame $e^\pm_{\varepsilon}=\exp\left\{{\pm
iu_\alpha^{\varepsilon}}\right\}e^\pm$ such that
$ D_A e^\pm_{\varepsilon} =\pm i\alpha d\theta_{\varepsilon}\otimes e_\pm
$.
 Normalize  by setting $u_\alpha^{\varepsilon}(1)=0$, $u^0_\alpha(z)\equiv 0$, and set
$\dot u_\alpha=\dot u_\alpha^\varepsilon\bigr|_{\varepsilon=0}$.  Then
$$ i\dot u_\alpha(z) = \frac{\alpha}{2}\left( \frac{\dot w}{z}-\overline{\frac{\dot
w}{z}}\right)\ .
$$
\end{Lem}

\begin{proof} Differentiate to obtain:
$$ D_A e^+_{\varepsilon}= idu_\alpha^{\varepsilon}\otimes e^+_{\varepsilon}+
i\alpha d\theta\otimes e^+_{\varepsilon}\
\Longrightarrow\ du_\alpha^{\varepsilon}= \alpha d\theta_{\varepsilon}-\alpha d\theta\ .
$$
 Hence, $d\dot u_\alpha=\alpha d\dot\theta$.  Now:
\begin{equation*}
d\theta_{\varepsilon}=\frac{1}{2i}\left(\frac{dw}{w}-\overline{\frac{dw}{w}}\right)
\ \Longrightarrow\ d\dot\theta =\frac{1}{2i}\left(d\left(\frac{\dot
w}{z}\right)-d\left(
\overline{\frac{\dot w}{z}}\right)\right)\ .
\end{equation*} 
Since $\dot w(1)=0$ by the normalization (\ref{E:normalization}), the result
follows upon integration.
\end{proof}

We now choose one more smooth cut-off function in addition to $\varphi_0$ and $\varphi_1$: 
let
$\varphi_2$ be smooth,  supported on $\Delta$ and is identically 1 on $\Delta_{1/3}$.  Set
$u_{\ab}^\varepsilon=u_\alpha^\varepsilon-u_\beta^\varepsilon$, and define the complex gauge
transformation:
$$
 g_{\ab}^{\varepsilon}=
          \left(    {\begin{matrix} \exp\left({i\varphi_2 u_{\ab}^{\varepsilon}+
\varphi_1(\alpha-\beta)\log|w_{\varepsilon}|}\right) & 0
\\
                         0 & \exp\left({-i\varphi_2 u_{\ab}^{\varepsilon}-
\varphi_1(\alpha-\beta)\log|w_{\varepsilon}|}\right)
                       \end{matrix}  }
          \right) \ .
$$
 Now consider the family of singular connections
$\gamma_{\varepsilon}=g_{\varepsilon} A_\alpha$, where the action of the complex gauge
transformation is with respect to the complex structure $w_{\varepsilon}$. 
The derivative in the $e^+$ direction, for example, is given by
 (see 
Thm.\ 
\ref{T:firstvariation}):
\begin{equation*}
\left(\dot\gamma^{0,1}\right)_{e^+}=
\frac{\partial}{\partial \bar z}\left[ i\varphi_2\dot
u_{\ab}+\varphi_1\frac{\alpha-\beta}{2}\left(\frac{\dot w}{z}+\overline{\frac{\dot w}{z}}
\right)\right]
 -2(\alpha-\beta)\tilde\nu \frac{\partial}{\partial z}\left(\varphi_1\log|z|^2\right)\ .
\end{equation*}
 By Lemma \ref{L:frametwisting} this is:
$$
\left(\dot\gamma^{0,1}\right)_{e^+}=
\frac{\partial}{\partial \bar z}\left[
\frac{\varphi_2(\alpha-\beta)}{2}\left\{ \left(\frac{\dot w}{z}-\overline{\frac{\dot w}{z}}
\right) +\varphi_1\left(\frac{\dot w}{z}+\overline{\frac{\dot w}{z}}
\right)\right\}\right]-2(\alpha-\beta)\tilde\nu \frac{\partial}{\partial
z}\left(\varphi_1\log|z|^2\right)\ .
$$
 By the choice of cut-off functions it is easily verified that the support of
$\dot\gamma^{0,1}$ lies in the annulus $1/6\leq |z|\leq 2/3$.  The continuity  follows
from this expression.
\end{proof}

%%%%%%%%%%%%%%%%%%%%%%%%%%%%%%%%%%%%%%%%%%%%%%%%%%%%%%%%%%

\subsection{Differentiability}  

%%%%%%%%%%%%%%%%%%%%%%%%%%%%%%%%%%%%%%%%%%%%%%%%%%%%%%%%%%

We combine Thm.\ \ref{T:firstvariation} and Prop.\ \ref{P:critical} to prove:

\begin{Thm} \label{T:diff}
Let $\sigma^\ast_{\varepsilon}$ be a continuously differentiable family of metrics representing a
path
$[\sigma^\ast_{\varepsilon}]\in\T(g,1)$.  Then for each $[A]\in \Ra$,
$\pi_\ab^{[\sigma^\ast_{\varepsilon}]}[A]$ is a continuously differentiable path in $\Rb$.
\end{Thm}

As an immediate consequence:

\begin{Cor}
Given $0<\beta<\alpha< 1/2$, the universal Hecke correspondence $\tilde\pi_\ab:\widetilde
\Ra\to\widetilde\Rb$ is continuously differentiable.  The same is true for $\beta=0$ on the
preimage of the irreducible representations.
\end{Cor}

\noindent \emph{Proof of Thm. \ref{T:diff}}.  First, notice that by choosing a rational number
$\beta<k/n<\alpha$ and using Lemma \ref{L:composition} it suffices to prove the result for
rational holonomies.  Second, by Def. \ref{D:foliation1}  and Prop.\ \ref{P:critical}  it suffices to
show that given a continuously differentiable path of connections $A_{\varepsilon}\in \A_\beta$, the
path
$[A_\varepsilon(\infty)]\in \Rb$ is also continuously differentiable.  Finally, by passing to a
branched cover as in Sec. \ref{S:covers}, it suffices to prove the result for closed surfaces.  We
continue with the notation as in Sec. \ref{S:covers}.

Let $\widehat A_{\varepsilon}\in \widehat\A$ be a continuously differentiable path. Note that by
Prop. \ref{P:branchedreducibles} and Lemma 
\ref{L:branchedflow} we may assume the $\widehat A_\varepsilon$ are stable. Since
$\widehat A_\varepsilon(\infty)=g_{\varepsilon}\widehat A_{\varepsilon}$, where $g_{\varepsilon}$ is
a complex gauge transformation, it suffices by the first variational formula (Thm.
\ref{T:firstvariation}) to show that $h_{\varepsilon}$ is a smooth family of complex gauge
transformations.  This can be achieved by the implicit function theorem as follows:  consider the
map
\begin{align*}
f:\Omega^0(\sqrt{-1}\ad_0\widehat E)&\times\widehat\A_{stable}\times\met \lra
\Omega^0(\sqrt{-1}\ad_0\widehat E)\\
f(u,\widehat A,\hat\sigma)&=\sqrt{-1}\ast_{\hat\sigma}\dbar_{\widehat
A}^{\hat\sigma}\left(e^{-u}\partial_{\widehat A}^{\hat\sigma} e^u\right)\ .
\end{align*}
It is easily verified that the map $f$ is smooth, and:
$$
\left(\delta_u f\right)_{(u,\widehat A,\hat\sigma)}(\delta u)=\Delta^{\hat\sigma}_{e^u\widehat
A}(\delta u)\ .
$$
By the
implicit function theorem, the solution $u=u(\widehat A,\hat\sigma)$  depends smoothly on $(\widehat
A,\hat \sigma)$, and this completes the proof.

%%%%%%%%%%%%%%%%%%%%%%%%%%%%%%%%%%%%%%%%%%%%%%%%%%%%%%%%%%%%%%%%%%%%%%

\section{Eigenvalue and Eigenfunction Estimates}    \label{S:eigenvalue}

%%%%%%%%%%%%%%%%%%%%%%%%%%%%%%%%%%%%%%%%%%%%%%%%%%%%%%%%%%%%%%%%%%%%%%%%

%%%%%%%%%%%%%%%%%%%%%%%%%%%%%%%%%%%%%%%%%%%%%%%%%%%%%%%%%%

\subsection{Estimates}   \label{S:estimates}

%%%%%%%%%%%%%%%%%%%%%%%%%%%%%%%%%%%%%%%%%%%%%%%%%%%%%%%%%%

We shall need  eigenfunction and eigenvalue estimates 
on sections of a  vector bundle  $V$ equipped
with a Riemannian metric, a metric connection $A$, and Laplace operator $\Delta_A$.  
  In the
following, we shall assume
$X$ is a compact surface, possibly with boundary, with a smooth metric of area ${\mathfrak a}$.
We define (cf.\   \cite[eq.\   
(0.4)]{Li}) the Sobolev constant ${\mathfrak s}_1$ to be the supremum over all constants $s$
satisfying:
$ s\inf_{a\in\RBbb}\ \Vert f-a\Vert_2^2\leq \Vert df\Vert_1^2
$, for all smooth functions $f$.  In case $\partial X\neq \emptyset$, we define another
constant associated to the Dirichlet problem: let ${\mathfrak s}_2$  be the supremum over
all constants $s$ satisfying:
$ s\,  \Vert f\Vert_2^2\leq \Vert df\Vert_1^2
$, for all smooth compactly supported functions $f$.  In this section we prove the following:

\begin{Thm}  \label{T:efbounds1} There is a universal constant $C$ with the following
property:  let $\varphi$ be an eigensection of $V$ with eigenvalue $\lambda$.  Then:
\begin{enumerate}
\item[(1)]  If $\partial X=\emptyset$ or if $\partial X\neq \emptyset$ and $\varphi$ satisfies
Neumann boundary conditions, then:
$$
\Vert \varphi\Vert^2_\infty\leq C\left(
{\mathfrak a}^{-1}+\left(\frac{4\lambda}{\mathfrak{s}_1}\right)^2{\mathfrak a}\right)
\Vert \varphi\Vert_2^2\ .
$$
\item[(2)]  If $\partial X\neq \emptyset$ and $\varphi$ satisfies Dirichlet boundary conditions,
then:
$$
\Vert \varphi\Vert^2_\infty\leq C\left(\frac{4\lambda}{\mathfrak{s}_2}\right)^2{\mathfrak a}\Vert
\varphi\Vert_2^2\ .
$$
\end{enumerate}
\end{Thm}

We also need a lower bound on the growth of eigenvalues. While heat kernel estimates as in
\cite{CL} might perhaps give more precise estimates, we shall only need the following:

\begin{Thm} \label{T:evbounds1} There is a universal constant $C$ with the following
property: 
\begin{enumerate}
\item[(1)]  If $\lambda_k$ denotes the $k$-th eigenvalue for sections of $V$ for the closed or
Neumann boundary problem, then:
$$ k\leq C\rk V\left( 1+\frac{4\lambda_k {\mathfrak a}}{\mathfrak{s}_1}\right)^3\ .
$$
\item[(2)]  If $\mu_k$ denotes the $k$-th eigenvalue for sections of $V$ for the  Dirichlet
boundary problem, then:
$$ k\leq C\rk V\left(
\frac{4\mu_k {\mathfrak a}}{\mathfrak{s}_2}\right)^3\ .
$$
\end{enumerate}
\end{Thm}

Thm.'s \ref{T:efbounds1} and \ref{T:evbounds1} stated above are 
generalizations of the results of P. Li  \cite{Li} (see also \cite{CL}) on eigenvalues and
eigenfunctions of the Laplacian on forms.
 For the sake of
completeness, however, we shall sketch the important steps involved in the proofs of the results
for sections. We point out that the main difference in  the estimates is the
extra term ${\mathfrak a}^{-1}$ in Thm.\  \ref{T:efbounds1}, Part (1), and the $1$ in Thm.\ 
\ref{T:evbounds1}, Part (1).  These are basically due to the lack of an \` a priori lower bound on
the eigenvalues in terms of Sobolev constants. The key estimate  is the following  (cf.\  
\cite[eq.\    (2.7) and (3.1)]{Li}):

\begin{Lem}  \label{L:keyestimate} Let $\varphi$ be a smooth section of $V$ and $A$ a
smooth metric connection on $V$.  If
$\partial X\neq \emptyset$, we assume that $\varphi$ satisfies either Dirichlet or Neumann
boundary conditions.  Then for any
$\alpha > 1$,
$$
\int_X
|\varphi|^{2\alpha-2}\langle\varphi,\Delta_A\varphi\rangle\geq\frac{2\alpha-1}{\alpha^2}
\int_X\left|d|\varphi|^\alpha\right|^2\ .
$$
\end{Lem}

\begin{proof} Using $ 2\langle\varphi, D_A\varphi\rangle=d|\varphi|^2
$, we have for $\alpha>1$:
$$
\left\langle D_A\left(|\varphi|^{2\alpha-2}\varphi\right), D_A\varphi\right\rangle
=\frac{\alpha-1}{2}|\varphi|^{2\alpha-4}\langle d|\varphi|^2,
d|\varphi|^2\rangle+|\varphi|^{2\alpha-2} |D_A\varphi|^2\ .
$$ By Kato's inequality, the right hand side above is
$$
\geq \frac{\alpha-1}{2}|\varphi|^{2\alpha-4}\langle d|\varphi|^2,
d|\varphi|^2\rangle+|\varphi|^{2\alpha-2}\left| d|\varphi|\right|^2
=\frac{2\alpha-1}{\alpha^2}\left| d|\varphi|^\alpha\right|^2\ .
$$ The lemma then follows from integration by parts.
\end{proof}

\noindent \emph{Proof of Thm.\  \ref{T:efbounds1}}.  We shall only prove Part (1), the proofs
of the other statements being similar.  Assume
$\Delta_A\varphi=\lambda\varphi$.  Applying (\ref{E:df}) to $f=|\varphi|^\alpha$ and using Lemma
\ref{L:keyestimate}, we obtain:
$$
\lambda\int_X |\varphi|^{2\alpha}\geq \frac{2\alpha-1}{\alpha^2}\frac{\mathfrak{s}_1}{4}
\left(
{\mathfrak a}^{-1/2}\left(\int_X|\varphi|^{4\alpha}\right)^{1/2}-{\mathfrak a}^{-1}\int_X|
\varphi|^{2\alpha}\right)\ ,
$$ or,
$$
\left( {\mathfrak a}^{-1/2}+\frac{\alpha^2}{2\alpha-1}\frac{4\lambda
{\mathfrak a}^{1/2}}{\mathfrak{s}_1}\right)^{1/2\alpha}\Vert\varphi\Vert_{2\alpha}\geq
\Vert\varphi\Vert_{4\alpha}\ .
$$ Set $\alpha=2^k$ to obtain, by iteration,
$$
\Vert\varphi\Vert_2\prod_{j=0}^k\left( {\mathfrak a}^{-1/2}+\frac{2^{2j}}{2^{j+1}-1}\frac{4\lambda
{\mathfrak a}^{1/2}}{\mathfrak{s}_1}\right)^{1/2^{j+1}}\geq
\Vert\varphi\Vert_{2^{k+2}}\ .
$$ Letting $k\to\infty$ yields
$\  \displaystyle
\Vert\varphi\Vert^2_\infty\leq \Vert\varphi\Vert_2^2 {\mathfrak a}^{-1}\prod_{j=0}^\infty\left(1
+\frac{2^{2j}}{2^{j+1}-1}\frac{4\lambda {\mathfrak a}}{\mathfrak{s}_1}\right)^{1/2^j}
$. The result then follows from the following simple:
\begin{Lem}  \label{L:orange} Fix $\gamma >0$.  Then for $\beta>1$ there is a constant
$C(\beta)$ independent of $\gamma$ such that:
$$
\prod_{j=0}^\infty\left(1 +\frac{\gamma \beta^{2j}}{2\beta^{j}-1}\right)^{1/\beta^j}\leq
C(\beta)\left(1+\gamma^{\beta/\beta-1}\right)\ .
$$
\end{Lem}

\noindent \emph{Proof of Thm.\  \ref{T:evbounds1}}.  Again, we shall concentrate on Part
(1).  Let the first $k$ eigenvalues and eigensections be denoted
$0\leq \lambda_1\leq\ldots\leq \lambda_k$, and $\varphi_1,\ldots,\varphi_k$, respectively.

\begin{Lem}  \label{L:mango} For any $\psi\in\text{span}\, \{\varphi_1,\ldots,\varphi_k\}$
and any $l\geq 1$,
$$
\frac{ \Vert\psi\Vert_{2^{l+2}}  }{ {\mathfrak a}^{1/2^{l+2}}  }
\leq
\left( 1+\frac{4\lambda_k {\mathfrak a}}{\mathfrak{s}_1}\frac{2^{2l}}{2^{l+1}-1}\right)^{1/2^{l+1}-1}
\frac{ \Vert\psi\Vert_{2^{l+1}}  }{ {\mathfrak a}^{1/2^{l+1}}  }\ .
$$
\end{Lem}

\begin{proof} By Lemma \ref{L:keyestimate} for $\alpha=2^l$ and  (\ref{E:df}) we have:

\begin{align*}
\frac{(2^{l+1} -1)}{2^{2l}}\frac{\mathfrak{s}_1}{4}\biggl\{
{\mathfrak a}^{-1/2}&\Vert\psi\Vert_{2^{l+2}}^{2^{l+1}} - {\mathfrak a}^{-1} \Vert\psi\Vert_{2^{l+1}}^{2^{l+1}}
\biggr\} 
\leq
\int_X |\psi|^{2^{l+1}-2}\langle\psi,\Delta_A\psi\rangle   \\ &\leq
\left\{ \int_X |\psi|^{2^{l+1}}\right\}^{(2^{l+1}-1)/2^{l+1}}
\left\{ \int_X |\Delta_A\psi|^{2^{l+1}}\right\}^{1/2^{l+1}}\ .
\end{align*}

\noindent   Write $\psi=\sum_{j=1}^k c_j\varphi_j$. Now it follows as in \cite[Lemma
17]{Li} that there is a subset $J\subset\{1,\ldots,k\}$ such that:
\begin{align*}
\left\{ \int_X |\Delta_A\psi|^{2^{l+1}}\right\}^{1/2^{l+1}} &\leq
\left\{ \int_X |\Delta_A\psi|^{2^{l+2}}\right\}^{1/2^{l+2}}{\mathfrak a}^{1/2^{l+2}}  \leq
  \left\{ \int_X \bigl|\sum_{j=1}^k \lambda_j
c_j\varphi_j\bigr|^{2^{l+2}}\right\}^{1/2^{l+2}}  {\mathfrak a}^{1/2^{l+2}}  
\\ &\leq   \lambda_k   \left\{ \int_X \bigl|\sum_{j\in J} 
c_j\varphi_j\bigr|^{2^{l+2}}\right\}^{1/2^{l+2}}  {\mathfrak a}^{1/2^{l+2}}  \leq
\lambda_k\Vert\psi\Vert_{2^{l+2}}  {\mathfrak a}^{1/2^{l+2}}   \ .
\end{align*} Therefore,
$$ {\mathfrak a}^{-1/2}\left( \frac{  \Vert\psi\Vert_{2^{l+2}}  }{  \Vert\psi\Vert_{2^{l+1}} }
\right)^{2^{l+1}-1}
\leq \frac{4\lambda_k}{\mathfrak{s}_1}\frac{2^{2l}}{2^{l+1}-1}{\mathfrak a}^{1/2^{l+2}}+ {\mathfrak a}^{-1}
\frac{  \Vert\psi\Vert_{2^{l+1}}  }{ \Vert\psi\Vert_{2^{l+2}}   }\ .
$$ Using the fact that $\Vert\psi\Vert_{2^{l+1}}\leq {\mathfrak a}^{1/2^{l+2}}\Vert\psi\Vert_{2^{l+2}}$,
we have:
$$
\left( \frac{\Vert\psi\Vert_{2^{l+2}}}{ \Vert\psi\Vert_{2^{l+1}}} \right)^{2^{l+1}-1}
\leq \left( 
{\mathfrak a}^{-1/2}+\frac{4\lambda_k}{\mathfrak{s}_1}\frac{2^{2l}}{2^{l+1}-1}{\mathfrak a}^{1/2}\right)
{\mathfrak a}^{1/2^{l+2}}\ ,
$$ from which the lemma follows.
\end{proof}

\noindent  For each $l\geq 1$, choose $\psi_l\in\text{span}\,
\{\varphi_1,\ldots,\varphi_k\}$ so that 
$
 \Vert \varphi\Vert_{ 2^{l+2}} \Vert\psi_l\Vert_2 \leq 
 \Vert \psi_l\Vert_{ 2^{l+2}}\Vert\varphi\Vert_2
$ for all $\varphi\in\text{span}\, \{\varphi_1,\ldots,\varphi_k\}$.  Then repeated
application of Lemma \ref{L:mango} gives:
$$
\frac{ \Vert \varphi\Vert_\infty }{\Vert\varphi\Vert_2}
\leq
\frac{ \Vert \psi_0\Vert_4 }{\Vert\psi_0\Vert_2}
\frac{1}{{\mathfrak a}^{1/4}}
\prod_{l=1}^\infty\left( 1+\frac{4\lambda_k {\mathfrak a}}{\mathfrak{s}_1}
\frac{2^{2l}}{(2^{l+1}-1)}\right)^{1/2^{l+1}-1}  \ .
$$ From Lemma \ref{L:keyestimate}  and  (\ref{E:df}) we also have:
$$
\frac{ \Vert \psi_0\Vert_4^2 }{\Vert\psi_0\Vert_2^2}
\leq {\mathfrak a}^{-1/2}
\left( 1+\frac{4\lambda_k {\mathfrak a}}{\mathfrak{s}_1}\right)\ .
$$ Using Lemma \ref{L:orange} and the fact that $1/(2^{l+1}-1)\leq 1/2^l$ we conclude that
for all
$\varphi\in\text{span}\, \{\varphi_1,\ldots,\varphi_k\}$
$$
\frac{ \Vert \varphi\Vert_\infty^2 }{\Vert\varphi\Vert_2^2}
\leq C  {\mathfrak a}^{-1}
\left( 1+\frac{4\lambda_k {\mathfrak a}}{\mathfrak{s}_1} \right)^3\ .
$$ Part (1) of Thm.\  \ref{T:evbounds1} now follows from the following:
\begin{Lem}[{\cite[Lemma 11]{Li}}] Let $W$ be a finite dimensional subspace of $L^2(V)\cap
C^0$.  Then there is $\varphi\in W$ such that
$
\ \dim W\, \Vert \varphi\Vert_2^2\leq (\rk V)\, {\mathfrak a} \Vert\varphi\Vert_\infty^2\ .
$
\end{Lem}

%%%%%%%%%%%%%%%%%%%%%%%%%%%%%%%%%%%%%%%%%%%%%%%%%%%%%%%%%%

\subsection{Conic Degeneration}  \label{S:conicdegeneration}

%%%%%%%%%%%%%%%%%%%%%%%%%%%%%%%%%%%%%%%%%%%%%%%%%%%%%%%%%%

Recall the notion of a  cone metric on a manifold $C(Y)=(0,1)\times Y$ (cf.\   \cite{Ch1},
\cite{Ch2},
\cite{JW1}):  this is a metric of the form $ds^2=dr^2+r^2 \tilde \sigma$, where $\tilde
\sigma$ is a (smooth) metric on $Y$.  An $n$-dimensional manifold $X$ with metric $\sigma$
on
$X\setminus\{p\}$ is said to have a cone metric if for some choice of $(Y,\tilde \sigma)$
and some neighborhood $U$ of
$p$ in $X$,
$U\setminus\{p\}$ is isometric to $C(Y)$; in this case we call $p$ a cone point.  We
generalize this notion to that of a cone double point, by which we mean locally the union
of two copies of
$C(Y)$ with the singularity identified.  For example, for surfaces it is natural to view
such a singularity as arising from the following family of degenerating metrics:  on the
cylinder
$ C=\left\{ (x,y) : -1\leq x\leq 1,\ 0\leq y\leq 2\pi\right\}\bigr/ \left\{ (x,0)\sim
(x,2\pi)\right\}
$, define, for $0\leq \ell\leq 1$, the  metrics $ds_\ell^2 = dx^2 +
\left(\ell+(1-\ell)x^2\right)\kappa^2dy^2$. The parameter $0<\kappa\leq 1$, introduced here for
convenience,  is a fixed cone angle.  The metric for
$\ell=0$ has a cone double point as described above.  We shall refer to a family of metrics
$\sigma(\ell)$ on a  compact, connected surface
$X$ as a \emph{conic degenerating family} if
$X$ contains a finite collection of disjoint  cylinders $C_1,\ldots, C_k$, all of   which
are isometric to
$C(\ell)=(C, ds_\ell^2)$, and
$\sigma(\ell)$ converges smoothly to a Riemannian metric on $X\setminus C_1\cup\cdots\cup C_k$.

Eigenvalue problems are still well-behaved on compact manifolds
with cone metrics.  We shall be considering the following situation:
 let $A$ be a connection on a (real) vector bundle $V\to X$ with Riemannian metric. If $X$
has cone points, then we allow $A$ to be singular at those points. It will suffice to
assume further
 that with respect to some choice of orthonormal  frame and conformal coordinates
$(r,\theta)$ near the cone point, $A$ has the form $b\otimes d\theta $, where $b$ is
a constant diagonal matrix. Consider the operator $\Delta_A=D_A^\ast D_A$ acting on sections
$\varphi$ of
$V$ satisfying
$\varphi,\ D_A \varphi,\ \Delta_A\varphi$ in $L^2$.  Then by a generalization of  the
result in
\cite{Ch1}, Stokes theorem holds for sections of $V$ and the cone metric.  Hence, the $L^2$
extension of $\Delta_A$ is self-adjoint.  Moreover, the spectrum of $\Delta_A$ is discrete
with finite multiplicity, and the eigenfunctions are smooth away from the singularities (cf.\  
\cite{Ch2}).

We will be using three important properties of conic degenerating families.  The first follows
from a direct computation:

\begin{Prop}  \label{P:curvaturebound} Let $(X,\sigma(\ell))$ be a conic degenerating family
of surfaces.  Then the eigenvalues of the Ricci curvature tensor are
uniformly bounded by a constant times $1/\ell$.
\end{Prop}

\begin{proof}  We briefly sketch the computation:  it clearly suffices to
compute the Ricci tensor
${R_k}^i$ for $C(\ell)=(C, ds_\ell^2)$.  We use coordinates $x_1=x$, $x_2=y$, so that:
$$
\sigma_{ij}=\left(\begin{matrix} 1 & 0 \\ 0 & \left(\ell+(1-\ell)x^2\right)\kappa^2 
\end{matrix}\right)\ .
$$
Then one finds for the Christoffel symbols: $\Gamma^1_{11}=\Gamma^2_{11}=
\Gamma^1_{12}=\Gamma^2_{22}=0$, and
$$
\Gamma^1_{22}=-\frac{1}{2}\frac{\partial\sigma_{22}}{\partial x}=-(1-\ell)x\kappa^2 \quad ,
\quad
\Gamma^2_{12}=\frac{1}{2}\sigma^{22}\frac{\partial\sigma_{22}}{\partial x}=
\frac{(1-\ell)x}{\ell+(1-\ell)x^2}\ .
$$
The operator
$$
{R_k}^i=\sigma^{lj}\left(\frac{\partial \Gamma^i_{jl} }{\partial x^k}
-\frac{\partial \Gamma^i_{jk} }{\partial
x^l}+\Gamma^\mu_{jl}\Gamma^i_{\mu k}-\Gamma^\mu_{jk}\Gamma^i_{\mu l}
\right)\ .
$$
From this, one verifies ${R_1}^2={R_2}^1=0$, and 
\begin{align*}
{R_1}^1&=\sigma^{22}\left(\frac{\partial \Gamma^1_{22} }{\partial x}
-\Gamma^2_{12}\Gamma^1_{22}
\right)=\frac{-\ell(1-\ell)}{(\ell+(1-\ell)x^2)^2}\ ,\\
{R_2}^2&=-\frac{\partial \Gamma^2_{12} }{\partial x}
-\Gamma^2_{12}\Gamma^2_{12}
=\frac{-\ell(1-\ell)}{(\ell+(1-\ell)x^2)^2}\ .
\end{align*}
The result follows.
\end{proof}

 The second fact is a comparison between conic degeneration and the well-known
``plumbing" construction used to study holomorphic degenerating families of Riemann
surfaces (cf.\  
\cite{DW2}).  Consider the annuli given by:
$ \varUpsilon_\varepsilon=\left\{ (z,w)\in \CBbb^2 : |z|,|w|\leq 1, \;
zw=\varepsilon\right\}
$. This is a holomorphic family for $\varepsilon$ in a punctured disk, but we will take
$\varepsilon$ to be real.  The central fiber $\varepsilon=0$ is a ``pinched" annulus -- two
disks with coordinates $z$ and $w$ identified at $z=w=0$.  We shall be interested in
metrics on
$\varUpsilon_\varepsilon$ which are conformal with respect to the complex structure induced
from
$\CBbb^2$ and which degenerate to cone metrics
$ds_0^2$ on
$\varUpsilon_0$ proportional to metrics of the following form:  for some
$0<\kappa\leq 1$ and in the coordinate $z$,

\begin{equation}  \label{E:conemetric} 
ds_0^2=\kappa^2|z|^{2(\kappa-1)}|dz|^2
\end{equation}

\noindent (cf.\   \cite[Lemma 6.1]{JW1}).  These can always be constructed:

\begin{Prop}    \label{P:conicdegeneration} Fix $\kappa$, $0<\kappa\leq 1$.  Then there is
a function $\varepsilon=\varepsilon(\ell)$ depending on $\kappa$ and  a conformal conic
degenerating family of metrics on
$\varUpsilon_{\varepsilon(\ell)}$ converging to the metric 
(\ref{E:conemetric}).  Furthermore, the parameters
$\varepsilon$ and $\ell$ are related by the following bound:
$$
\frac{\ell^{2/\kappa}}{4^{1/\kappa}}  \leq \varepsilon(\ell)\leq \ell^{1/\kappa}\ ,
$$ for $\ell\leq 3/4$.
\end{Prop}

\begin{proof}
We solve for conformal coordinates $z=re^{i\theta}$ on the portion of
$\varUpsilon_{\varepsilon(\ell)}$ where $0\leq x\leq 1$.  Take $\theta=\theta(y)=y$, and
$r=f(x,\ell)$ for $f$ increasing, $f(1,\ell)=1$.  For such a solution, we may write
$ds^2_\ell=\sigma(z,\ell)|dz|^2$, or:
$$
\frac{dr^2}{(f')^2(x,\ell)}+(\ell+(1-\ell)x^2)\kappa^2
d\theta^2=\sigma(z,\ell)(dr^2+r^2d\theta^2)\ .
$$
where $f'(x,\ell)$ is the partial derivative with respect to $x$.  This implies:
$$
\frac{\partial}{\partial x}\log f(x,\ell)=\frac{1}{\kappa(\ell+(1-\ell)x^2)^{1/2}}\quad ,\quad
\sigma(z,\ell)=\frac{1}{(f')^2(x,\ell)}\ .
$$
Applying the initial condition, we find the solution
\begin{align*}
f(x,\ell)&=
\left[ \frac{(1-\ell)^{1/2}x}{1+(1-\ell)^{1/2}}+\frac{(1-\ell)^{1/2}}{1+(1-\ell)^{1/2}}
\sqrt{x^2+\frac{\ell}{1-\ell}}
\right]^{1/\kappa(1-\ell)^{1/2}} \\
\varepsilon(\ell)&= f^2(0,\ell)=\left[
\frac{\ell}{(1+(1-\ell)^{1/2})^2}\right]^{1/\kappa(1-\ell)^{1/2}}
\end{align*}
Notice that the relationship between $\varepsilon$ and $f$ in the second line follows from the fact
that the middle of the cylinder corresponds to
$|z|=|w|=\sqrt\varepsilon$. The convergence and bounds on $\varepsilon(\ell)$ easily follow.
\end{proof}

Lastly, we comment on the behavior of the Sobolev constants  under conic degeneration.  We quote
the following result:

\begin{Prop}[{\cite[Thm.\  2.4 and Prop.\  2.6]{JW1}}]   \label{P:cylindersobolevbound}
 For the conic degenerating
cylinder
$C(\ell)$ in Sec.\  \ref{S:conicdegeneration}, there is a constant $c$ such that
${\mathfrak s}_2(C(\ell))\geq c>0$ for all $\ell>0$.
\end{Prop}

This result bears on the inclusion $L^2_1\hookrightarrow L^4$.  
 Recall that for a smooth function $f$ the definition of
the Sobolev constant implies (see \cite[Lemma 4]{Li}):
\begin{equation}  \label{E:df}
\Vert df\Vert_2^2 \geq \frac{\mathfrak{s}_1}{4}\left( {\mathfrak a}^{-1/2}\Vert f\Vert_4^2-{\mathfrak a}^{-1}\Vert
f\Vert_2^2\right)\ ,
\end{equation}
where $\mathfrak a$ denotes the area of $X$.
For a conic degenerating family, $\mathfrak{s}_1\to 0$ if and only if there is a separating
pinching cylinder (see \cite[Cor.\  2.9]{JW1}).  However, by the proposition above, the
Sobolev constants $\mathfrak{s}_2$ of the component regions, i.e. the degenerating cylinders
and their complements, all remain bounded away from zero.  Thus, applying a cut-off
function to (\ref{E:df}), one easily proves

\begin{Prop}  \label{P:sobolevbound}
For a conic degenerating family $(X,\sigma(\ell))$ there is a constant $c>0$ independent of
$\ell$ such that for all smooth functions $f$ on $X$ and all $\ell>0$,
$$
\Vert f\Vert_{4;\sigma(\ell)}^2\leq c\,\bigl(\Vert df\Vert_{2;\sigma(\ell)}^2+\Vert
f\Vert_{2;\sigma(\ell)}^2\bigr)\ .
$$
\end{Prop}

%%%%%%%%%%%%%%%%%%%%%%%%%%%%%%%%%%%%%%%%%%%%%%%%%%%%%%%%%%%%%%%%%%%%%%%%%%%%%%%%%%%%%%%%%%

\subsection{Uniform Bounds}  \label{S:uniformbounds}

%%%%%%%%%%%%%%%%%%%%%%%%%%%%%%%%%%%%%%%%%%%%%%%%%%%%%%%%%%%%%%%%%%%%%%%%%%%%%%%%%%%%%%%%%%

If $(X,\sigma(\ell))$ is a conic degeneration, then a vector bundle $V\to X$ with connection
$A$ such that $A$ has the standard form $b\otimes d\theta$, $b$ constant diagonal, on the
cylinder
$C$, naturally defines a family of operators $\Delta_A$ on $(X,\sigma(\ell))$.  By the comments at
the end of Sec.\ \ref{S:conicdegeneration}, it is reasonable to ask whether  the eigenvalues and
eigenfunctions of
$\Delta_A$ on
$(X,\sigma(\ell))$  converge as $\ell\to 0$ to those of the limiting cone metric.  This is what we
call \emph{spectral convergence}.
The estimates from Thm.'s \ref{T:efbounds1} and \ref{T:evbounds1} are the key elements
needed to prove spectral convergence for this degenerating family (cf.\   \cite{JW1, JW2}).  The
first step is uniform $C^0$ bounds on eigenfunctions and uniform growth of eigenvalues:

\begin{Cor}  \label{C:efbounds2} 
If $(X,\sigma(\ell))$ is a  conic degenerating family with
vector bundle $V$ and connection $A$, then there are constants
$C_1$,
$C_2$ independent of $\ell$ such that if $\varphi(\ell)$ is a normalized eigensection of $V$ with
eigenvalue
$\lambda(\ell)$ then
$
\Vert\varphi(\ell)\Vert_\infty\leq  C_1+C_2\lambda^5(\ell)
$ for all $\ell$.
\end{Cor}

\begin{Cor}  \label{C:evbounds2} If $(X,\sigma(\ell))$ is a  conic degenerating family,
then there is a constant $C$ and an integer $N$, both  independent of $\ell$, such that if
$\lambda_k(\ell)$ denotes the $k$-th eigenvalue for the closed problem for sections of $V$,
then
$
\lambda_k(\ell)\geq C k^{1/3}
$ for all $\ell$ and all $k\geq N$.
\end{Cor}

 \noindent The
arguments in the references mentioned above then apply to give:

\begin{Thm}   For a  conic degenerating family we have spectral convergence for sections of
$V$.
\end{Thm}

\noindent Applying this to the particular case of a  unitary connection $A$ on $E$ and
the eigenvalue
$\lambda_1$ of
$\Delta_A$ acting on sections of $(\ad E)_0$, we have:

\begin{Cor}   \label{C:lambdabound}
Let $A$ be an irreducible  connection on $X$.  Then for a  conic
degenerating family
$(X,\sigma(\ell))$, $\lambda_1(\ell)\to 0$ if and only if $A$ is accidentally reducible.
\end{Cor}

\begin{proof}[Proof of Cor.\  \ref{C:efbounds2}] In the non-separating case, it follows
from \cite[Cor.\  2.9 and Thm.\  2.4]{JW1} that the Sobolev constant $\mathfrak{s}_1$
is bounded away from zero.  The result then follows from Part (1) of Thm.\ 
\ref{T:efbounds1}.  
In the separating case, following \cite{JW2} we  divide the degenerating surface into three
regions
$X^\pm$ and
$C$ with the  induced metrics from $\sigma(\ell)$ (denoted also by $\sigma(\ell)$).   Let
$\widehat C\subset C$ be a fixed subcylinder and $\widehat X^\pm$ the corresponding
complementary regions. By definition, the degenerating family $(\widehat C, \sigma(\ell))$
is isometric to the standard degenerating cylinder; therefore, by
Prop.\  \ref{P:cylindersobolevbound} it follows that the Sobolev constant
$\mathfrak{s}_2$ for
$(\widehat C, \sigma(\ell))$ is bounded away from zero.  The same is true for $(\widehat
X^\pm,
\sigma(\ell))$, since these form a smooth family of metrics on the surfaces with boundary
for
$\ell\geq 0$.

Choose a smooth cut-off function
$\eta$ on $X$ such that $\eta\equiv 1$ on $\widehat C$ and $\eta\equiv 0$ on $X^\pm$.
Suppose
$\varphi(\ell)$ is a normalized eigensection on $(X,\sigma(\ell))$ with eigenvalue
$\lambda(\ell)$. We denote by
$((1-\eta)\varphi(\ell))^\pm$ the restriction of $((1-\eta)\varphi(\ell))$ to $\widehat
X^\pm$. It clearly suffices to find bounds on $((1-\eta)\varphi(\ell))^\pm$ and
$\eta\varphi$ separately. Notice that $\Delta_A^m((1-\eta)\varphi(\ell))^\pm$ and
$\Delta_A^m\eta\varphi$  have
$L^2$ bounds depending on $m$ and $\lambda(\ell)$, but otherwise independent of $\ell$. 
More precisely, the bound may be taken of the form $C_1+C_2\lambda^m(\ell)$.  This is
because
$d\eta$ and $d(1-\eta)$ are supported in $\widehat X^\pm$ where the higher derivatives of
$\varphi$ may be uniformly bounded by an application of the elliptic estimate.

Let $\psi^\pm_k(\ell)$ and $\psi^C_k(\ell)$ denote normalized eigensections with eigenvalues
$\mu^\pm_k(\ell)$ and $\mu^C_k(\ell)$ for the Dirichlet problems on $\widehat X^\pm$ and
$\widehat C$, respectively.  Consider the Fourier expansions:
$$ ((1-\eta)\varphi(\ell))^\pm=\sum_{k=1}^\infty a_k^\pm(\ell)\psi^\pm_k(\ell)
\quad ,\quad
\eta\varphi(\ell)=\sum_{k=1}^\infty b_k(\ell)\psi^C_k(\ell)\ .
$$ Since the Sobolev constants for the individual pieces are bounded away from zero, it
follows from Thm.\  \ref{T:efbounds1} that  we have uniform bounds:
$$
\Vert \psi^\pm_k(\ell) \Vert_\infty \leq C \mu^\pm_k(\ell) \quad , \qquad
\Vert \psi^C_k(\ell) \Vert_\infty \leq C \mu^C_k(\ell)
\; .
$$ Here the constant $C$ may be chosen independent of $k$ and $\ell$. Hence, it suffices to
show that there is some $B$ of the required form satisfying:
\begin{align}
\label{first}
\sum_{k=1}^\infty \left|  a_k^\pm(\ell) \right| \mu^\pm_k(\ell) &\leq B \ ,\\
\label{second}
\sum_{k=1}^\infty \left| b_k(\ell)\right|
\mu^C_k(\ell)&\leq B\ ,
\end{align} for all $\ell$.  Consider (\ref{second}).  By definition, we have:
\begin{align*} b_k(\ell) =
\int_{C} \langle\eta\varphi(\ell), \psi^C_k(\ell)\rangle d\sigma(\ell)
&=\frac{1}{(\mu^C_k(\ell))^m}
\int_{C} \langle\eta\varphi(\ell),\Delta^m_A
\psi^C_k(\ell)\rangle d\sigma(\ell) \\ &=\frac{1}{(\mu^C_k(\ell))^m}
\int_{C} \langle\Delta^m_A\left(\eta\varphi(\ell)\right),
\psi^C_k(\ell)\rangle d\sigma(\ell)\ .
\end{align*} Since we assume uniform $L^2$ bounds on
$\Delta^m_A\left(\eta\varphi(\ell)\right)$, we obtain a uniform bound
$$
\mu^C_k(\ell) \left| b_k(\ell)\right|
\leq \left(C_1+C_2\lambda^m(\ell)\right)(\mu^C_k(\ell))^{1-m}\ ,
$$ where the constants $C_1, C_2$ may be chosen independent of $\ell$. By Part (2) of Thm.\ 
\ref{T:evbounds1} above we have
$\mu^C_k(\ell)\geq C' k^{1/3}$ for some constant $C'$ independent of $k$ and
$\ell$. Hence,
$$
\mu^C_k(\ell) \left| b_k(\ell)\right|
\leq \left(C^\prime_1+C^\prime_2\lambda^m(\ell)\right) k^{(1-m)/3}\ ,
$$ where the constants $C_1^\prime, C_2^\prime$ may be chosen independent of $k$ and 
$\ell$. Since  the sum over $k$ of the terms on the right hand side converges for $m\geq
5$, the desired bound in (\ref{second}) is obtained.  The proof for (\ref{first}) is
similar.
\end{proof}

\begin{proof}[Proof of Cor.\  \ref{C:evbounds2}] In the non-separating case, it follows
from \cite[Cor.\  2.9 and Thm.\  2.4]{JW1} that the Sobolev constant $\mathfrak{s}_1$
is bounded away from zero.  The result then follows from Part (1) of Thm.\ 
\ref{T:evbounds1}.
In the separating case, consider the three regions $X^\pm$, $C$ as in the proof of Cor.\ 
\ref{C:evbounds2} above.  
 By domain monotonicity, it suffices to prove the result for the Neumann spectra for
$(X^\pm, \sigma(\ell))$ and $(C,\sigma(\ell))$.   Since
$(X^\pm,
\sigma(\ell))$ form a smooth family of metrics on the surfaces with boundary, it follows
that the Sobolev constants
$\mathfrak{s}_1$ for these regions are uniformly bounded away from zero; hence, by Part (1) of
Thm.\ 
\ref{T:evbounds1} the Neumann spectra of these surfaces has uniform growth as in the
statement of the corollary. As above, the Sobolev constant
$\mathfrak{s}_2$ for $(C,\sigma(\ell))$ is also uniformly bounded away from zero.  Hence,
by Part (2) of Thm.\ 
\ref{T:evbounds1} the Dirichlet spectra of these surfaces has uniform growth as in the
statement of the corollary.   On the other hand, because of the form of the connection on
$C$, the Neumann spectrum  may be bounded below by the Dirichlet
spectrum after shifting the index by two as in [W].  This completes the proof.
\end{proof}

%%%%%%%%%%%%%%%%%%%%%%%%%%%%%%%%%%%%%%%%%%%%%%%%%%%%%%%%%%%%%%%%%%%%%%%%%%%%%

\subsection{Heat Kernel Estimates}  \label{S:heatkernel}

%%%%%%%%%%%%%%%%%%%%%%%%%%%%%%%%%%%%%%%%%%%%%%%%%%%%%%%%%%%%%%%%%%%%%%%%%%%%%

It is well-known that eigenvalue and eigenfunction estimates produce estimates on
solutions to the linear heat equation.  This will also be useful in the non-linear Yang-Mills
flow.  The result we need is the following:

\begin{Thm}   \label{T:heatkernel}
Let $(X,\sigma(\ell))$ is a  conic degenerating family and fix $T>0$.  There is a constant $C$
depending on $T$ but independent of
$\ell$ such that if $v(t,x)$ is the solution to the heat equation with initial conditions
$v(0,x)=v_0$, then $\sup_x|v(t,x)|\leq C\Vert v_0\Vert_{2;\sigma(\ell)}$ for all $t\geq T$.
\end{Thm}

\begin{proof}   The solution $v(t,x)$ may be written:
$$ v(t,x)=\sum_{i=0}^\infty e^{-t\lambda_i} a_i\varphi_i(x)\ ,
$$ 
where $|a_i|\leq \Vert v_0\Vert_{2;\sigma(\ell)}$. By Cor.\  \ref{C:efbounds2} for the
ordinary Laplacian on functions,  there are constants
$C_1$ and
$C_2$ independent of $\ell$ such that:
$$
\sup_x|v(t,x)|\leq \Vert v_0\Vert_{2;\sigma(\ell)}\sum_{i=0}^\infty
e^{-t\lambda_i(\ell)}\left(C_1+C_2\lambda_i^5(\ell)\right)\ .
$$ 
Thus, we must show that for fixed $T_0>0$ there is a constant $ C$ depending only on
$T_0$, $C_1$, and $C_2$ such  that for all $t\geq T_0$ and all $\ell>0$, 
$$
\sum_{i=0}^\infty e^{-t\lambda_i}\left(C_1+C_2\lambda_i^5(\ell)\right)\leq \widetilde C\ .
$$
 First choose $\Lambda_0$ such that for all $\Lambda\geq \Lambda_0$ we
have
$\Lambda^6(C_1+C_2\Lambda_i^5)e^{-\Lambda T}\leq 1$.  Now by Cor.\ 
\ref{C:evbounds2} we can find a constant
$C_3$ such that $\lambda_{i}\geq C_3 i^{1/3}$, for $i$ sufficiently large, say $i\geq N$,
independent of $\ell$.  We further prescribe $N$ such that for
$i\geq N$,
$\lambda_i\geq
\Lambda_0$.  Finally, let:
$$ C_4=\sup_{\Lambda>0}(C_1+C_2\Lambda_i^5)e^{-\Lambda T}\ .
$$ Then for $t\geq T$,
\begin{align*}
\sum_{i=0}^\infty e^{-t\lambda_i}(C_1+C_2\Lambda_i^5) & =
\sum_{i=0}^{N-1} e^{-t\lambda_i}(C_1+C_2\Lambda_i^5)  +\sum_{i=N}^\infty
e^{-t\lambda_i}(C_1+C_2\Lambda_i^5) \\
&\leq N C_4 + \frac{1}{C_3^6}\sum_{i=N}^\infty
\frac{1}{i^2}\ .
\end{align*}
 This proves the result.
\end{proof}

%%%%%%%%%%%%%%%%%%%%%%%%%%%%%%%%%%%%%%%%%%%%%%%%%%%%%%%%%%%%%%%%%%%%%%%%%%%%%

\section{Proof of the Main Theorem}  \label{S:proof}

%%%%%%%%%%%%%%%%%%%%%%%%%%%%%%%%%%%%%%%%%%%%%%%%%%%%%%%%%%%%%%%%%%%%%%%%%%%%%

%%%%%%%%%%%%%%%%%%%%%%%%%%%%%%%%%%%%%%%%%%%%%%%%%%%%%%%%%%

\subsection{Outline of the Proof}  \label{S:outline}

%%%%%%%%%%%%%%%%%%%%%%%%%%%%%%%%%%%%%%%%%%%%%%%%%%%%%%%%%%

Let $[\sigma^\ast(\ell)]$ be a degeneration in $\T_{aug.}(g,1)$ to a nodal Riemann surface
with conformal structure $[\sigma^\ast(0)]$ associated to a collection $\Phi$ of simple closed
curves. Recall that $\R^\Phi\subset\R$  denotes the $\Phi$-accidentally reducible representations. In
this section we are going to show that for a given
$[A]\in\Ra\setminus \R^\Phi$, and $|\beta-\alpha|<\varepsilon_0$, where
$\varepsilon_0$ is sufficiently small  as in Sec.\  \ref{S:foliation},
 
\begin{equation}   \label{E:piconvergence}
\lim_{\tau\to\infty} \pi_\ab^{[\sigma^\ast(\ell)]}[A]=\pi_\ab^{[\sigma^\ast(0)]}[A]\ .
\end{equation}

By standard compactness arguments, this will suffice to prove the Main Theorem.
As representatives for the degenerating conformal structure we may choose lifts
$\sigma^\ast(\ell)\to\sigma^\ast(0)$ to be a conic degeneration as in Sec.\ 
\ref{S:conicdegeneration}  with cone angle $0<\kappa<1$. The fact that $\kappa$ may be chosen 
strictly less than $1$ will be important (see (\ref{E:choice})).   Given $[A]\in\Ra\setminus\R^\Phi$,
let
$g_\ab(A)$ denote the twists of $A$ in the standard form around the point $p$.  Let $A(\ell,t)$,
$A(0,t)$ denote the Yang-Mills flow of $g_\ab(A)=A_0$ with respect to $\sigma^\ast(\ell)$,
$\sigma^\ast(0)$, respectively.
Since a conjugacy class of flat connections is completely determined by its holonomy around
all homotopy classes of closed curves on $X^\ast$, our strategy of proof will be to first
show that the holonomies of $A(\ell,\infty)$ converge for curves supported away from the
pinching cylinders. This does not suffice, however.  Indeed, this statement, combined with the
results on Simpson's flow from Sec.\ 
\ref{S:gaugetheory}, show only that the limiting holonomies of $A(\ell,\infty)$ around the
pinching cylinders return to the initial holonomies of $A_0$ as $\ell\to 0$.  This would
still allow for the possibility of a change of framing, or gluing parameters, across the
cylinder (cf.\ the discussion preceding Def.\ \ref{D:infinitefoliation}).  So the second part of the
proof is to show that the holonomies
\emph{across} the pinching cylinders, as measured with respect to the framing coming from $A_0$, are
very nearly trivial. We present these two results as Thm.'s \ref{T:holonomy} and
\ref{T:transverseholonomy} below:

\begin{Thm}  \label{T:holonomy}
For any set $\{\Xi_j\}_{j=1}^N$ of  closed curves supported in  $\pi_1(X_0^\ast)$ we have
$$
\left\{\hol_{\Xi_j} A(\ell,\infty)\right\}_{j=1}^N\lra
\left\{\hol_{\Xi_j} A(0,\infty)\right\}_{j=1}^N\ ,
$$
modulo overall conjugation by $SU(2)$.
\end{Thm}

Next, recall from Thm.\  \ref{T:lift}  that the manner by which a connection
$[A(0,\infty)]$ produces a point in $\Rb$ is to use the initial framings.  Consider a
cylinder $(C, ds_\ell^2)$ in $X^\ast$ on which the twisted initial connection $A_0$ is flat.  We
may choose a unitary frame $\{e_1,e_2\}$ such that in terms of the coordinates $(x,y)$ in Sec.\ 
\ref{S:conicdegeneration}, $A_0$ has the form
$
d_{A_0}=d+\diag(i\gamma, -i\gamma) dy
$.
We fix this frame once and for all throughout the degeneration.  Now for a small transverse
arc $\Gamma_\varepsilon=\{ (x,y)\in C : -\varepsilon\leq x\leq \varepsilon,\ y=y_0\}$ and flat
connection
$A(\ell,\infty)$, we measure the holonomy $\hol_{\Gamma_\varepsilon}\left(
A(\ell,\infty)\right)$ by parallel translating the frame $\{e_1,e_2\}$ along
$\Gamma_\varepsilon$.  For example, notice that by our choice of lift $[A(0,\infty)]$ in
Thm.\ 
\ref{T:lift}, $\hol_{\Gamma_\varepsilon}\left( A(0,\infty)\right)=\I$ for any choice of
$\varepsilon$.

Since any closed curve on $X^\ast$ may be written, up to homotopy, as a concatenation of
curves of the form $\Xi$ in Thm.\  \ref{T:holonomy} and transverse arcs
$\Gamma^i_{\varepsilon}$, one for each component  $c_i\in\Phi$, 
we see that (\ref{E:piconvergence}) will follow from Thm.\  \ref{T:holonomy} and the
following:

\begin{Thm}     \label{T:transverseholonomy}
For any $\delta>0$ there is  $\varepsilon_0 >0$ and $\ell_0$ such that for all
$\Gamma^i_{\varepsilon_0} $ and all $\ell\geq \ell_0$,
$$
\left| \hol_{\Gamma^i_{\varepsilon_0} }\left(A(\ell,\infty)\right)-\I\right| < \delta\ .
$$
\end{Thm}

\noindent  The proof of Thm.\  \ref{T:holonomy} will occupy the next two subsections, and the
proof of Thm.\  \ref{T:transverseholonomy} will be given in subsection
\ref{S:transverseholonomy}.  One of the key ingredients in the proofs is the $C^0$ estimate for
the metric and the curvature found in Cor.\  \ref{C:conicexpdecay} below.

%%%%%%%%%%%%%%%%%%%%%%%%%%%%%%%%%%%%%%%%%%%%%%%%%%%%%%%%%%

\subsection{Proof of Theorem \ref{T:holonomy}}      \label{S:holonomy}

%%%%%%%%%%%%%%%%%%%%%%%%%%%%%%%%%%%%%%%%%%%%%%%%%%%%%%%%%%

Throughout this section, $\Omega$ is an open set with
compact closure in $X^\ast\setminus{\Phi}$.  As before, we fix $[A]\in\Ra$.  Choose a
representative $A$, and let $A_0=g_\ab(A)$ be the twist of $A$ at $p$.  We assume that
$\beta = k/n$, where $k,n$ are positive coprime integers, and $n$ is odd.

We denote by $A(\ell,t)$ the solution of the Yang-Mills flow
(\ref{E:ymflow})-(\ref{E:initialcondition}),   by $h(\ell,t)$ the solution to the
non-linear heat equation (\ref{E:metricflow})-(\ref{E:metricinitialcondition}) on
$(X^\ast, \sigma^\ast(\ell))$, and by
$h(0,t)$  the solution to the same equations on
$X^\ast_0=X^\ast\setminus{\Phi}$ with the degenerate metric $\sigma^\ast(0)$, the holomorphic
structure on the bundle being determined by $A_0$.

As a first step, we show that as we degenerate the metric on $X^\ast$, $h(\ell,t)$ converges
to
$h(0,t)$ uniformly on compact sets.  In the following, $p$ is any number strictly greater than
1.  Also, $\ast_\ell$ and $\ast_0$ will denote the Hodge $\ast$'s on $X^\ast$ and
$X^\ast\setminus\Phi$ with respect to the metrics $\sigma^\ast(\ell)$ and $\sigma^\ast(0)$,
respectively.

\begin{Prop} \label{P:loghbound}
Given  $T>0$, 
$
 \log h(\ell, T)\to\log h(0,T)
$ weakly in $L^p_{2,loc.}$.  In particular, the convergence is strong in $C^1(\Omega)$.
\end{Prop}

\begin{proof}  In Sec.\  \ref{S:C1bound}, we will obtain $C^0$ bounds for $h(\ell,t)$ and
$\ast_\ell F_{A(\ell,t)}$ independent of $\ell$ (see Cor.\  \ref{C:conicexpdecay}).  Assuming
these results, since 
$$
h^{-1/2}(\ell,t)F_{A(\ell,t)}h^{1/2}(\ell,t)=F_{A_0}+\dbar_{A_0}\left(h^{-1}(\ell,t)\partial_{A_0}
h(\ell,t)\right)\ ,
$$
standard elliptic estimates imply $L^p_{2,loc.}$ estimates for $h(\ell,t)$ uniform in $\ell$ and
$0\leq t\leq T$.  By eq.\    (\ref{E:metricflow}), this implies an $L^p_{2,1,loc.}$ estimate on
$X^\ast\times[0,T]$, uniform in $\ell$, where the $1$ refers to the time derivative.  It follows
that $h(\ell, t)$ converges to some $\tilde h(0,t)$ weakly in $L^p_{2,loc.}$.  The uniform $C^0$
bounds imply that $\tilde h(0,t)$ and $\ast_0 F_{\dbar_{A_0},\tilde h(0,t)H_0}$ are also bounded
uniformly for $t\in [0,T]$.  The uniqueness part of Thm.\  \ref{T:simpson} shows that
$\tilde h(0,T)=h(0,T)$ as desired.  
\end{proof}

\begin{Cor}  \label{C:loghbound}
Given $\varepsilon >0$ there exists $T>0$ and $\ell_0=\ell_0(\varepsilon,T)>0$ such that for
$\ell\geq
\ell_0$,
$$
\Vert \log h(\ell, T)-\log h(0,\infty)\Vert_{C^1(\Omega)} < \varepsilon\ ,
$$
and similarly,
$$
\Vert  h^{1/2}(\ell, T)A_0- h^{1/2}(0,\infty)A_0\Vert_{C^0(\Omega)} < \varepsilon\ .
$$
\end{Cor}

\begin{proof}
By Thm.\   \ref{T:simpson},  $h(0,t)\to h(0,\infty)$ uniformly in $C^1(\Omega)$, so we can
choose
$T_1$ such that for all $t\geq T_1$,
$$
\Vert \log h(0,t)-\log h(0,\infty)\Vert_{C^1(\Omega)} < \varepsilon/2\ .
$$
Take $T\geq T_1$.  By Prop.\  \ref{P:loghbound}, $\exists\ \ell_0=\ell_0(\varepsilon,T)$ such
that for all $\ell\geq \ell_0$, 
$$
\Vert \log h(\ell,t)-\log h(0,T)\Vert_{C^1(\Omega)} < \varepsilon/2\ .
$$
The result follows.
\end{proof}

Now let $A(\ell,t)$ be the solution of the Yang-Mills flow on $X^\ast$ as before.  Write
\begin{equation}
A(\ell,\infty)=\tilde g(\ell,t)A(\ell,t)\quad ,\quad
\tilde h(\ell,t)=\tilde g^\ast(\ell,t)\tilde g(\ell,t)\ .
\end{equation}

\begin{Prop}   \label{P:C1bound}
Given $\varepsilon>0$, there is a $T_0>0$ independent of $\ell$ such that for all $t\geq
T_0$, $|\log\tilde h(\ell,t)|_{C^1(\Omega)} <\varepsilon$.
\end{Prop}

We shall also prove this result in the following subsection.  Here we show how Prop.\ 
\ref{P:C1bound} implies Thm.\  \ref{T:holonomy}.
Let $T=\max(T_0,T_1)$, where $T_0$ and $T_1$ are as in Prop.\  \ref{P:C1bound} and
Prop.\  \ref{P:loghbound}, respectively.  Let $\ell_0$ be chosen as in Cor.\ 
\ref{C:loghbound}, and choose
$\ell\geq
\ell_0$.  Notice that both $h^{1/2}(\ell,T)A_0$ and $\tilde h^{1/2}(\ell,T)A(\ell,\infty)$
are real gauge equivalent to $A(\ell,T)$.  We therefore can write:
$$
\tilde h(\ell,T)A(\ell,\infty)=k(\ell,T)h^{1/2}(\ell,T)A_0\ ,
$$
where $k(\ell,T)$ is a real gauge transformation.  By Cor.\  \ref{C:loghbound}:
$$
\left\Vert \tilde
h^{1/2}(\ell,T)A(\ell,\infty)-k(\ell,T)h^{1/2}(0,\infty)A_0
\right\Vert_{C^0(\Omega)}\leq
c\, \varepsilon\ ,
$$
for $c$ depending only on the Sobolev embedding $L^p_1(\Omega)\hookrightarrow
L^\infty(\Omega)$, and may be taken independent of $\ell$.
On the other hand, by Prop.\  \ref{P:C1bound}:
$$
\left\Vert \tilde
h^{1/2}(\ell,T)A(\ell,\infty)-A(\ell,\infty)\right\Vert_{L^\infty(\Omega)}\leq
\varepsilon\ .
$$
It follows that:
$$
\left\Vert
A(\ell,\infty)-k(\ell,T)h^{1/2}(0,\infty)A_0\right\Vert_{L^\infty(\Omega)}\leq
(c+1)\varepsilon\ .
$$
Since $A(\ell,\infty)$ and $k(\ell,T)h^{1/2}(0,\infty)A_0$ are $C^0$-close, their
holonomies around the $\Xi_j$ are also close.  Finally, since
$k(\ell,T)h^{1/2}(0,\infty)A_0$ and $A(0,\infty)$ are real gauge equivalent, the
theorem  follows.

%%%%%%%%%%%%%%%%%%%%%%%%%%%%%%%%%%%%%%%%%%%%%%%%%%%%%%%%%%

\subsection{Proof of Proposition \ref{P:C1bound}}   \label{S:C1bound}

%%%%%%%%%%%%%%%%%%%%%%%%%%%%%%%%%%%%%%%%%%%%%%%%%%%%%%%%%%

We begin with some preliminary results:

\begin{Prop}  \label{P:closetoflat}
Let $[A]\in\Ra$ be as before.  Then given $\varepsilon_0>0$ there is $\delta >0$ such that
for $|\beta-\alpha| <\delta$ there exists a twist $A_0 =g_\ab(A)$ and a flat connection
$A_\infty$  with the standard form of holonomy $\beta$ around $p$, and such that
$
\Vert A_0-A_\infty \Vert_{4;\sigma(\ell)} <\varepsilon_0
$.
\end{Prop}

\begin{proof}
Let $C$ denote a component of the pinching region and $A_0$ the twist of $A$ at $p\not\in C$.  Let
$h$ denote the flow at $\infty$  for eq.'s (\ref{E:metricflow})-(\ref{E:metricinitialcondition})
associated to
$\dbar_{A_0}$, the initial hermitian structure $H_0$, and the degenerate metric $\sigma(0)$
on $X^\ast_0$.  Let $\widetilde A_\infty$ denote the hermitian connection associated to
$\dbar_{A_0}$ and $hH_0$.  We may assume that $A_0$ is in the standard form
$d+\diag(i\gamma,-i\gamma)d\theta$ on a slightly larger cylinder $C_1\supset\overline C$.  Let
$C_2$ be an open cylinder such that $\overline C\subset C_2\subset\overline C_2\subset C_1$. 
Since the initial curvature $\Vert\ast F_{A_0}\Vert_\infty$ can be made arbitrarily small for
$|\beta-\alpha|$ small, and since $A_0$ is stable on $X_0^\ast$, it follows (cf.\   \cite{Si1}) that
$\Vert\log h\Vert_{C^1(V)}$, and hence also $\Vert A_0-\widetilde A_\infty\Vert_{C^0(V)}$, can be
made arbitrarily small for a relatively compact set $V\subset X^\ast_0$.  By \cite[Lemma 2.7]{DW1}
there is a real gauge transformation $g$ such that $g(\widetilde A_\infty)=A_0$ on $C_2$. 
Moreover, it is clear from the proof of that lemma that by taking $V$ so that $X\setminus V\subset
C$ and using the fact that
$\Vert A_0-\widetilde A_\infty\Vert_{C^0(V)}$ is small, we may conclude that $\Vert\log
g\Vert_{C^0(C_1\setminus C_2)}$ is small.  By bootstrapping we find that $\Vert\log
g\Vert_{L^4_1(C_1\setminus C_2)}$ is small; hence, we can extend $g$ to a real gauge
transformation of $E$ over $X^\ast$ with $g\equiv \I$ on $X\setminus C_1$ and $\Vert
g-\I\Vert_{L^4_1(X\setminus C_2)}$ small.  Set $A_\infty=g(\widetilde A_\infty)$.  Then $A_\infty$
extends to a connection over $X^\ast$, and the desired estimate holds.
\end{proof}

In order to get our estimate for the metric in $X^\ast$, we need again to pass to branched
covers.  Assume that $p\in X^\ast\setminus X$ is outside the pinching region, and let
$q:\widehat X\to X$ be a regular cyclic branched cover of degree $n$, chosen so that all branch
points lie outside the pinching region.  We choose metrics $\hat\sigma(\ell)$ on $\widehat
X$ so that $\hat\sigma(\ell)\to\hat\sigma(0)$ is a conic degeneration, and 
with respect to the induced conformal structures from $\hat\sigma(\ell)$ and $\sigma(\ell)$,
the map $q$ is holomorphic.  Recall the map $\hat q:\A_{k/n}\to\widehat\A:A\mapsto\widehat
A$.

\begin{Lem}    \label{L:lambdabound}
Let $A\in\A_{k/n}$ be a flat connection which is lies outside the $\Phi$-accidental
reducibles.  Then there is a $\lambda$  such that
$\lambda_1\bigl(\Delta_{\widehat A}^{\hat\sigma(\ell)}\bigr)\geq \lambda >0$ for all
$\ell$.
\end{Lem}

\begin{proof}
According to Cor.\  \ref{C:lambdabound}, $\lambda_1\bigl(\Delta_{\widehat
A}^{\hat\sigma(\ell)}\bigr)\to 0$ if and only if $\widehat A$  is accidentally reducible on
$\widehat X_0$.  But this is ruled out by the assumptions and Prop.\ 
\ref{P:branchedreducibles}.
\end{proof}

The next step is to
reduce to the case of a closed surface. First, notice that it suffices to get $C^0$ estimates for
$\hat h(\ell,t)$ and $\ast_\ell F_{\widehat A(\ell,t)}$.  Also, by Lemma
\ref{L:branchedestimate}, and shrinking
$\Omega$ slightly to avoid branch points, it suffices to show:

\begin{equation}     \label{E:C1bound}
\bigl\Vert \log\widetilde{\hat h}(\ell,t)\bigr\Vert_{C^1(\widehat\Omega)} <\varepsilon\ .
\end{equation}

Note also that by Prop.\  \ref{P:closetoflat} we have a flat connection $\widehat A_\infty$
such that
$\Vert
\widehat A_\infty-\widehat A_0\Vert_4 <\varepsilon_0$, where $\varepsilon_0$ can be taken
arbitrarily small.  Furthermore, by Lemma \ref{L:lambdabound} we may assume that 
 $\lambda_1\bigl(\Delta_{\widehat A}^{\hat\sigma(\ell)}\bigr)\geq \lambda >0$ uniformly in
$\ell$.  It follows by the fundamental estimate \cite[Prop.\  7.2 -- second method]{R}
that if $\varepsilon_0$ is less than a universal constant $\varepsilon_1$, depending only on
$\lambda$, then
\begin{equation}  \label{E:initiallowerbound}
\left\Vert D_{\widehat A_0}^{\ast_\ell} F_{\widehat A_0}\right\Vert_{2,\hat\sigma(\ell)}\geq
c\left\Vert\ast_\ell F_{\widehat A_0}\right\Vert_{2,\hat\sigma(\ell)}\ .
\end{equation}
Since from now on we will work on a closed surface $\widehat X$, we henceforth omit the hats
from the notation.  Again, $\sigma(\ell)\to\sigma(0)$ is a conic degeneration of a
closed surface
$X$.  We begin with the following:

\begin{Prop}   \label{P:conicexpdecay}
Let $A(\ell,t)$ be a solution of the Yang-Mills flow with initial condition $A_0$, $\Vert\ast_\ell
F_{A_0}\Vert_\infty<B$, and  such that the
fundamental estimate
$$
\left\Vert D_{ A(\ell,t)}^{\ast_\ell} F_{ A(\ell,t)}\right\Vert_{2;\sigma(\ell)}\geq
c\,\left\Vert\ast_\ell F_{ A(\ell,t)}\right\Vert_{2;\sigma(\ell)}\ ,
$$
holds for a constant $c$ independent of $\ell$ and for $0\leq t <\widehat T\leq\infty$.  Then
there are constants $c_1$, $c_2$ depending on $c$ and $B$, but  independent of $\ell$ and
$\widehat T$, such that for all $0\leq
t\leq\widehat T$:
\begin{enumerate}
\item[(1)] 
$\left\Vert\ast_\ell F_{ A(\ell,t)}\right\Vert_{\infty}\leq c_1e^{-ct/2}$; and,
\item[(2)]  $\left\Vert h(\ell,t)\right\Vert_\infty\leq c_2$.
\end{enumerate}
\end{Prop}

\begin{proof} For (1), first note that from (\ref{E:ymgradient}):
\begin{equation}   \label{E:derivative}
\begin{aligned}
\frac{d}{dt}\left\Vert\ast_\ell F_{ A(\ell,t)}\right\Vert_{2;\sigma(\ell)}&=-\frac{1}{2}
\left\Vert D_{ A(\ell,t)}^{\ast_\ell} F_{ A(\ell,t)}\right\Vert^2_{2;\sigma(\ell)}
\left\Vert\ast_\ell F_{ A(\ell,t)}\right\Vert^{-1}_{2;\sigma(\ell)} \\
&\leq -\frac{c}{2}\left\Vert D_{ A(\ell,t)}^{\ast_\ell} F_{ A(\ell,t)}\right\Vert_{2;\sigma(\ell)}\
,
\end{aligned}
\end{equation}
from which we obtain $\left\Vert\ast_\ell F_{ A(\ell,t)}\right\Vert_{2;\sigma(\ell)}\leq
c_1e^{-ct/2}$.  Let $u(t,x)=\left| \ast_\ell F_{A(\ell,t)}\right|(x)$ be the pointwise
norm.  By
\cite[Prop.\  16]{Do2} it follows that $\dot u+\Delta_\sigma u\leq 0$, where
$\Delta_\sigma$ is the ordinary Laplacian with respect to the metric $\sigma$. By the maximum
principle, we can first choose $c_2$ so that (1) holds for $0\leq t\leq 1$.  Now we rechoose $c_2$
according to the heat kernel bounds of Thm. \ref{T:heatkernel} and again apply the maximum
principle.

For (2) consider the metric flow (\ref{E:metricflow}):
$$ h^{-1}\frac{dh}{dt}=-\sqrt{-1}\ast_\ell F_{hH_0}\quad , \quad h(0)=\I \ .$$
Multiplying the equation through by $h$ and taking traces, we find:
$$
\frac{d}{dt}\log\tr h\leq \widetilde C \left\Vert\ast_\ell F_{h(t)H_0}\right\Vert_{\infty}\ ,
$$ for some numerical constant $\widetilde C$. By (1), the right hand side is uniformly integrable
on $0\leq t\leq \widehat T$.
 Since $\det h(t)=1$, the result follows.
\end{proof}

\begin{Prop} \label{P:coniclowerbound}
 Let $A_0$ be as before and $\Vert\ast_\ell
F_{A_0}\Vert_\infty<\varepsilon_2$.  Then there exists a constant $c$ independent of $\ell$
such that:
\begin{equation}  \label{E:lowerbound}
\left\Vert D_{ A(\ell,t)}^{\ast_\ell} F_{ A(\ell,t)}\right\Vert_{2;\sigma(\ell)}\geq
c\left\Vert\ast_\ell F_{ A(\ell,t)}\right\Vert_{2;\sigma(\ell)}\ ,
\end{equation}
for all $t\geq 0$, where $A(\ell,t)$ is the solution of the Yang-Mills flow 
(\ref{E:ymflow})  with initial condition $A(\ell,0)=A_0$, and with respect to the metric
$\sigma(\ell)$.
\end{Prop}

\begin{proof}
By (\ref{E:initiallowerbound}), there is $c>0$ such that:
$$
\left\Vert D_{ A(\ell,0)}^{\ast_\ell} F_{ A(\ell,0)}\right\Vert_{2;\sigma(\ell)}\geq
c\left\Vert\ast_\ell F_{ A(\ell,0)}\right\Vert_{2;\sigma(\ell)}\ ,
$$
uniformly in $\ell$.  Let 
$
{\mathcal J}=\left\{ t\in [0,\infty) : \ \text{the estimate (\ref{E:lowerbound}) holds on} \
[0,t]\right\}
$.
Then ${\mathcal J}\neq\emptyset$.  Let $T=\sup{\mathcal J}$.  We assume $T<\infty$, and
derive a contradiction.

\medskip
\noindent {\bf Claim 1.}  For the constant $c$ appearing in (\ref{E:lowerbound}), $\Vert
A(\ell,t)-A(\ell,0)\Vert_{2;\sigma(\ell)}< c^{-1}\varepsilon_2$.

\bigskip\noindent   By  (\ref{E:derivative}),
$$
\frac{d}{dt}\left\Vert\ast_\ell F_{ A(\ell,t)}\right\Vert_{2;\sigma(\ell)}\leq -c\left\Vert D_{
A(\ell,t)}^{\ast_\ell} F_{ A(\ell,t)}\right\Vert_{2;\sigma(\ell)} =-c\left\Vert\frac{\partial
A}{\partial t}(\ell,t)\right\Vert_{2;\sigma(\ell)}\ ,
$$
from which the estimate in Claim 1 follows by integration.

\medskip
\noindent {\bf Claim 2.} Write $A(\ell,T)=g(\ell)A_0$, $h(\ell)=g(\ell)g^\ast(\ell)$.  Then
$\Vert g(\ell)\Vert_\infty$ and $\Vert g^{-1}(\ell)\Vert_\infty$ are uniformly bounded
independent of $\ell$. 

\bigskip\noindent   This is immediate from  Prop.\  \ref{P:conicexpdecay}.

\medskip
\noindent {\bf Claim 3.}  There exists a real gauge transformation $u(\ell)$ such that
$$
\Vert u(\ell)A(\ell,T)-A_0\Vert_{L^2_1(A_0);\sigma(\ell)} < c\, \varepsilon_2\ ,
$$
where $c$ is independent of $\ell$.

\bigskip\noindent
Write $A(\ell,T)=g(\ell)A_0$ as before.  By Claim 1,   $\Vert
g^{-1}(\ell)\bar\partial_{A_0} g(\ell)\Vert_{2;\sigma(\ell)}\leq c\, \varepsilon_2$, and by
Claim 2, $\Vert
\bar\partial_{A_0} g(\ell)\Vert_{2;\sigma(\ell)}\leq c\varepsilon_2$ (in this as in the
following, $c$ will denote a generic constant independent of $\ell$).  By elliptic
regularity and the fact that $\lambda_1\left(\Delta_{A_0}^{\sigma(\ell)}\right)\geq
\lambda>0$ (see Prop.\  \ref{P:branchedreducibles}) it follows that 
\begin{align*}
\Vert\partial_{A_0} g(\ell)\Vert^2_{2;\sigma(\ell)}
&\leq \Vert d_{A_0} g(\ell)\Vert^2_{2;\sigma(\ell)}
\leq \Vert g(\ell)^\perp\Vert^2_{L^2_1(A_0);\sigma(\ell)} \\
&\leq \lambda^{-1}\Vert \bar\partial_{A_0}g(\ell)^\perp\Vert^2_{2;\sigma(\ell)}
= \lambda^{-1}\Vert \bar\partial_{A_0}g(\ell)\Vert^2_{2;\sigma(\ell)}
\leq c\,\varepsilon_2\ ,
\end{align*}
where $g(\ell)^\perp$ is the $L^2$-orthogonal projection to the perp space of $\ker
d_{A_0}$.  This result, combined with Claim 2, implies that $\Vert
h^{-1}(\ell)\partial_{A_0}h(\ell)\Vert_{2;\sigma(\ell)}\leq c\,\varepsilon_2$.  On the other
hand,
$$
g^{-1}(\ell)F_{A(\ell,t)}g(\ell)=F_{A_0}+\bar\partial_{A_0}
\left(h^{-1}(\ell)\partial_{A_0}h(\ell)\right)\ .
$$
Hence, by again applying elliptic regularity,
$
\Vert h^{-1}(\ell)\bar\partial_{A_0}h(\ell)\Vert_{L^2_1(A_0);\sigma(\ell)}\leq
c\,\varepsilon_2$.
Since the connection defined from $\bar\partial_{A_0}$ and $h(\ell)H_0$ is real gauge
equivalent to $A(\ell,T)$, the proof of Claim 3 is complete.

For notational simplicity, set $\altA(\ell,0)=u(\ell)A(\ell,T)$.  Then:
\begin{equation} \label{E:tilde}
\left\Vert \altA(\ell,0)-A_0\right\Vert_{L^2_1(A_0);\sigma(\ell)}< c\, \varepsilon_2\ .
\end{equation}

\medskip
\noindent {\bf Claim 4.}  Let $\altA(\ell,t)$ denote the Yang-Mills flow with initial
condition $\altA(\ell,0)$ and with respect to the metric $\sigma(\ell)$.  Then there is
$\delta>0$ and a real gauge transformation $v(\ell,t)$ such that:
$$
\left\Vert v(\ell,t)\altA(\ell,t)-A_0\right\Vert_{4;\sigma(\ell)}<  \varepsilon_1/2\ ,
$$
for $0<t<\delta$ and all $\ell$.  Here, $\varepsilon_1$ is the universal constant so that
R\aa de's estimate (\ref{E:rade}) holds.

\bigskip\noindent
The proof of Claim 4 will be accomplished in three stages:

\noindent (i)  Given $\eta >0$, $\exists\ \delta_1 >0$ such that for $0\leq t\leq \delta_1$
and all $\ell$, $\Vert \log\alth(\ell,t)\Vert_\infty <\eta$, where $\alth(\ell,t)$ is the
solution of (\ref{E:metricflow})-(\ref{E:metricinitialcondition}) with respect to the
holomorphic structure defined by $\bar\partial_{\altA(\ell,0)}$. This follows by Prop.\ 
\ref{P:conicexpdecay}, as in Claim 2.

\noindent (ii)  Given $\eta >0$, $\exists\ \delta_2 >0$ such that for $0\leq t\leq \delta_2$
and all $\ell$, $\Vert \altA(\ell,t)-\altA(\ell,0)\Vert_{2;\sigma(\ell)}<\eta$. This obtains
from the following inequalities:
\begin{align*}
\Vert \altA(\ell,t)-\altA(\ell,0)\Vert_{2;\sigma(\ell)}&=
\left\Vert\int_0^t\frac{d}{dt}\,\altA(\ell,t)dt\right\Vert_{2;\sigma(\ell)}\leq
\int_0^t\left\Vert\frac{d}{dt}\,\altA(\ell,t)\right\Vert dt   \\
&=\int_0^t\left\Vert
D^{\ast_\ell}_{\altA(\ell,t)}F_{\altA(\ell,t)}\right\Vert_{2;\sigma(\ell)} dt
\leq
t^{1/2}\left(\int_0^t\left\Vert
D^{\ast_\ell}_{\altA(\ell,t)}F_{\altA(\ell,t)}\right\Vert^2_{2;\sigma(\ell)} dt
\right)^{1/2}\\
&=t^{1/2}\left(   \left\Vert\ast_\ell F_{\altA(\ell,0)}\right\Vert^2_{2;\sigma(\ell)}
-\left\Vert\ast_\ell F_{\altA(\ell,t)}\right\Vert^2_{2;\sigma(\ell)}\right) \leq c\,
t^{1/2}\varepsilon\ .
\end{align*}

\noindent (iii)  Let $\altA(\ell,t)=\altg(\ell,t)\altA(\ell,0)$,
$\alth(\ell,t)=\altg^\ast(\ell,t)\altg(\ell,t)$, and $\delta=\min(\delta_1,\delta_2)$.  As
in Claim 3, we first obtain
$\Vert\bar\partial_{\altA(\ell,0)}\altg(\ell,t)\Vert_{2;\sigma(\ell)}\leq c\,\eta$ for
$0\leq t\leq \delta$.  On the other hand, by choosing $\varepsilon_2$ sufficiently small,
it follows by (\ref{E:tilde}) that we may assume that
$\lambda_1\bigl(\Delta_{\altA(\ell,0)}^{\sigma(\ell)}\bigr)\geq
\lambda/2>0$.  As in Claim 3, we obtain for $0\leq t\leq \delta$,
$$
\Vert
\alth^{-1}(\ell,t)\bar\partial_{\altA(\ell,0)}\alth(\ell,t)\Vert_{L^2_1
\left({\altA(\ell,0)}\right);\sigma(\ell)}\leq c(\varepsilon_2+\eta)\ .
$$
By taking $\varepsilon_2$ and $\eta$ sufficiently small with respect to $\varepsilon_1$,
and applying this result together with  Kato's inequality and the uniform embedding
$L^2_1\hookrightarrow L^4$ for functions (see Prop.\  \ref{P:sobolevbound}), we obtain
Claim 4.

Now we are ready to complete the proof of the proposition.  By Claim 4, and by taking
$\varepsilon_0$ in Prop.\  \ref{P:closetoflat} sufficiently small, it follows that
$\Vert v(\ell,t)\altA(\ell,t)-A_\infty\Vert_{4;\sigma(\ell)}<\varepsilon_1$; hence,
$$
\left\Vert D_{ A(\ell,t)}^{\ast_\ell} F_{ A(\ell,t)}\right\Vert_{2;\sigma(\ell)}\geq
c\left\Vert\ast_\ell F_{ A(\ell,t)}\right\Vert_{2;\sigma(\ell)}\ ,
$$
for $0\leq t\leq T+\delta$, contradicting the assumption that $T=\sup{\mathcal J}$.
\end{proof}

The following corollary is an immediate consequence of Prop.'s \ref{P:conicexpdecay} and
\ref{P:coniclowerbound}, and completes the proof of the $C^0$ estimate:

\begin{Cor}  \label{C:conicexpdecay}
Let $A_0$ be as in Prop.\  \ref{P:coniclowerbound}.    Then
there are constants $c_1$, $c_2$, and $c_3$ independent of $\ell$ and $t$ such that
\begin{enumerate}
\item[(1)] 
$\left\Vert\ast_\ell F_{ A(\ell,t)}\right\Vert_{\infty}\leq c_1e^{-c_2t}$ for $0\leq
t\leq\widehat T$; 
\item[(2)]  $\left\Vert h(\ell,t)\right\Vert_\infty\leq c_3$.
\end{enumerate}
\end{Cor}

We now proceed with the proof of (\ref{E:C1bound}):

\begin{Cor}  
Let $\Omega\subset X$ be as before.  Write $A(\ell,t)=\tilde g(\ell,t)A(\ell,\infty)$,
$\tilde h(\ell,t)=\tilde g^\ast(\ell,t)\tilde g(\ell,t)$.  Given $\varepsilon >0$, $\exists
\ T_0>0$ independent of $\ell$ such that for all $t\geq T_0$, $\Vert \log\tilde
h(\ell,t)\Vert_{C^1(\Omega)} <\varepsilon$.
\end{Cor}

\begin{proof}
By Prop.'s  \ref{P:coniclowerbound} and  \ref{P:conicexpdecay},
we have $\left\Vert \ast_\ell
F_{A(\ell,t)}\right\Vert_\infty\to 0$ as $t\to\infty$ uniformly in $\ell$, and:
$$
\int_0^\infty\left\Vert \ast_\ell F_{A(\ell,t)}\right\Vert_\infty dt <\infty\ ,
$$
uniformly in $\ell$.    From this we
deduce as in Claims 1 and 2 that:
$$
\lim_{t\to\infty}\Vert A(\ell,t)-A(\ell,\infty)\Vert_{2;\sigma(\ell)}=0\ ,
$$
 and that $\Vert \tilde g(\ell,t)\Vert_{2;\sigma(\ell)}$ and 
$\Vert \tilde g^{-1}(\ell,t)\Vert_{2;\sigma(\ell)}$ are bounded, all uniformly in $\ell$. Hence,
$$\Vert \bar\partial_{A(\ell,\infty)}\tilde g(\ell,t)\Vert_{2;\sigma(\ell)}\ ,$$
can be made arbitrarily small independent of $\ell$.  Since $\lambda_1\bigl(
\Delta_{A(\ell,\infty)}^{\sigma(\ell)}\bigr)\geq
\lambda>0$ uniformly in $\ell$, it follows that $\Vert g(\ell,t)\Vert_{L^2_1;\sigma(\ell)}$
can be made arbitrarily small uniformly in $\ell$.  Finally, by using the curvature
estimate:
$$
\tilde g^{-1}(\ell,t)F_{A(\ell,t)}\tilde g(\ell,t)=\bar\partial_{A(\ell,\infty)}\left(
\tilde h^{-1}(\ell,t)\partial_{A(\ell,\infty)}\tilde h(\ell,t)\right)\ ,
$$
we obtain the corollary  by bootstrapping in $\Omega$.
\end{proof}

\noindent
This proves (\ref{E:C1bound}), and thus completes the proof of Prop.\ 
\ref{P:C1bound}.

%%%%%%%%%%%%%%%%%%%%%%%%%%%%%%%%%%%%%%%%%%%%%%%%%%%%%%%%%%%%%%%%%%%%%%%%%%%%%%%%%%%%%

\subsection{Proof of Theorem  \ref{T:transverseholonomy}}         \label{S:transverseholonomy}  

%%%%%%%%%%%%%%%%%%%%%%%%%%%%%%%%%%%%%%%%%%%%%%%%%%%%%%%%%%%%%%%%%%%%%%%%%%%%%%%%%%%%%

Let $h(\ell,\infty)$ denote the limit at infinite time of the solution to eq.'s
(\ref{E:metricflow})-(\ref{E:metricinitialcondition}).  Then up to real gauge,
$A(\ell,\infty)$ is $h^{1/2}(\ell,\infty)\cdot A_0$.  It suffices to obtain an estimate of the
form:
\begin{equation}  \label{E:estimateone}
\int_{\sqrt\varepsilon(\ell)}^{\varepsilon_0}\left|\imath_{\partial/\partial r}\left\{
h^{-1/2}(\ell,\infty)\partial_{A_0} h^{1/2}(\ell,\infty)-c.c.\right\}\right| dr\,
<\delta\ ,
\end{equation}
for sufficiently small $\varepsilon_0$ and $\ell$.  
In the above, the coordinate $r$ is $|z|$ for the conformal coordinates on $(C,ds_\ell^2)$ from
Prop.\  \ref{P:conicdegeneration}, and $\imath_{\partial/\partial r}$ denotes contraction of
the form in the $r$ direction.  The notation $c.c.$ means hermitian conjugation of the previous
term.   Also, note that we have taken the square root
$\sqrt\varepsilon$, since with respect to the conformal coordinates in Sec.\ 
\ref{S:conicdegeneration} this corresponds to a point in the middle of the cylinder, i.e. $x=0$. 
The estimate for the whole arc $\Gamma_{\varepsilon_0}^i$ follows from the estimates in the $z$
and $w$ coordinates separately. 

The square root on the gauge transformation $h$ is difficult to work with directly,
so we will eliminate it by using a $C^0$ bound.  Namely, to prove (\ref{E:estimateone}), is
suffices to prove:
\begin{enumerate}
\item[(1)]  There is $B$ such that $\Vert h(\ell,\infty)\Vert \leq B$ for all $\ell$; 
\item[(2)] $\displaystyle\int_{\sqrt\varepsilon(\ell)}^{\varepsilon_0}\left|
\imath_{\partial/\partial r}\left\{h^{-1}\partial_{A_0} h-c.c.\right\}\right| dr\, <\delta$
for $\varepsilon_0$ sufficiently small and $\sqrt\varepsilon(\ell)\leq \varepsilon$.
\end{enumerate}

\noindent    We further reduce this with the following:
\begin{Prop}  \label{P:estimatetwo}
Suppose that there is a constant $B$ such that (1) holds.  Suppose in addition that $\Vert
h^{-1}\partial_{A_0} h-c.c.\Vert_{2;\sigma(\ell)}\leq B$ for all $\ell$.  Then (2) also holds.
\end{Prop}

\begin{proof} Set $w= h^{-1}\partial_{A_0} h-c.c.$.  The first step of the proof is to show that
the hypotheses imply an estimate of the form:
\begin{equation}  \label{E:L4bound}
\Vert
w\Vert^4_{4;\sigma(\ell)}\leq B/\ell^2 \ ,
\end{equation}
where $B$ is independent of $\ell$.  Consider the function $|w|$.  By
Prop.\  \ref{P:sobolevbound}, it suffices to estimate $d|w|$.  By
Kato's inequality,
$
\left|d|w|\right|^2\leq \left|\nabla_{A_0}w\right|^2
$.
So it suffices to estimate the right hand side.  We now apply the Weitzenb\"ock formula for a
1-form $\omega$ with values in the self-adjoint bundle:
$$
\Delta_{A_0}\omega=-\nabla_{A_0}^\ast\nabla_{A_0}\omega+\left\{R,\omega\right\}+
\left\{F_{A_0},\omega\right\}\ .
$$
Now $F_{A_0}$ is uniformly bounded in any norm, say $C^1$, and we have
$\Delta_{A_0}h^{-1}\partial_{A_0}h=-\bar\partial_{A_0}^\ast F_{A_0}$; so $\Delta_{A_0} w$ is
bounded as well.  Therefore, the desired estimate is obtained, provided we can estimate the term
$\left\{R,w\right\}$.  The explicit formula can found in \cite{Wu}, and it
involves the operator ${R_k}^i=\sigma^{lj} {R_{klj}}^i$ (see \cite[p. 953]{Wu}).  Then
(\ref{E:L4bound}) follows from Prop.\  \ref{P:curvaturebound}.

Now we consider the problem on a cylinder $(C,ds_{\ell}^2)$.  Choose a local holomorphic frame
$\{f_i\}$ for $A_0$ adapted to $\{e_i\}$ as in Sec.\  \ref{S:ymflow}.   In conformal coordinates,
we may write
$$
h^{-1}\partial_{A_0}h=\sum_{n\in\ZBbb} c_n^{ij} z^n f_i\otimes f_j^\ast\otimes dz\ ,
$$
since $\bar\partial_{A_0}\left(h^{-1}\partial_{A_0}h\right)=0$.   The coefficients $c_n^{ij}$ depend
on
$\ell$, and we wish to estimate them.  For convenience, we set
$|c_{-1}|^2=|c_{-1}^{11}|^2+|c_{-1}^{22}|^2$.
  We assume the frame has been chosen such that $|f_i\otimes
f_j^\ast|=|z|^{\gamma_i-\gamma_j}$, where
$\gamma_1=\gamma$,
$\gamma_2=-\gamma$.  We will assume that $\gamma\neq 0$,  the argument in the case $\gamma=0$
being similar. Also, since we assume that $\sigma(\ell)$ converges to a cone metric of the type
(\ref{E:conemetric}), the integrals below will be carried out with respect to this metric.  A
simple computation shows that these estimates are valid.
 We have:
$$
\left|h^{-1}\partial_{A_0}h\right|^2=\sum_{m,n\in\ZBbb} c_m^{ij}\bar c_n^{ij} z^m\bar z^n
|z|^{\gamma_i-\gamma_j}|dz|^2\ .
$$
The $L^2$ bound on $w$ implies one on $|h^{-1}\partial_{A_0}h|$, which in turn implies that there is
a constant
$B$ independent of
$\ell$ such that:
$$
|c_m^{ij}|^2\int_{\varepsilon(\ell)\leq|z|\leq 1}|z|^{2m+\gamma_i-\gamma_j}|dz|^2\leq B\ .
$$
Therefore, there is another constant independent of $\ell$, which we also denote by $B$, such that:
\begin{enumerate}
\item[(i)]  For $m\geq 0$, $|c_m^{ij}|\leq B$.  Also, $|c_{-1}^{12}|\leq B$;
\item[(ii)]  For $m\leq -2$, or $(m,i,j)=(-1,2,1)$, $|c_m^{ij}|\leq B
\varepsilon(t)^{-m-\frac{1}{2}(\gamma_i-\gamma_j)-1}$;
\item[(iii)]  $|c_{-1}|^2\log(1/\varepsilon(\ell))\leq B$.
\end{enumerate}
We now apply this to:
$$
\int_{\varepsilon^{1/2}}^{\varepsilon_0}\left|\imath_{\partial/\partial
r}\left\{h^{-1}\partial_{A_0}h\right\}\right| dr
\leq \sum_{(m,i,j)\neq (-1,i,i)}|c_m^{ij}|
\int_{\varepsilon^{1/2}}^{\varepsilon_0}r^{m+\frac{1}{2}(\gamma_i-\gamma_j)} dr
+2|c_{-1}|\log(1/\varepsilon)
\ .
$$
By the estimates (i) and (ii) above, the first term on the right hand side may be made arbitarily
small for small $\varepsilon_0$, independent of $\ell$.  To estimate the second term, we use the
$L^4$ bound (\ref{E:L4bound}).  Because of the log term, it suffices to show that $|c_{-1}|$ vanishes
as some power of $\varepsilon$.  The estimate (iii) is not sufficient.  

To use the $L^4$ bound, we first isolate the $m=-1$ terms. Again, (\ref{E:L4bound}) implies a
similar bound on the $L^4$ norm of $h^{-1}\partial_{A_0}h|$.  Write:
$$
|h^{-1}\partial_{A_0}h|^2=\frac{|c_{-1}|^2+|c_{-1}^{12}||z|^{2\gamma}+
|c_{-1}^{21}||z|^{-2\gamma}}{|z|^2}
+ g\ ,
$$
Integrating over a subcylinder $C_a$, we have:
\begin{equation}     \label{E:firstcoefficient}
\int_{C_a}\frac{|c_{-1}|^4}{|z|^4}ds_0^2\leq B\int_{C_a} \left(|h^{-1}\partial_{A_0}h|^4
+g^2 +  
\frac{|c_{-1}^{12}|^4|z|^{4\gamma}}{|z|^4}+\frac{|c_{-1}^{21}|^4|z|^{-4\gamma}}{|z|^4}          
\right) ds_0^2\ .
\end{equation}
The subcylinder is given by $C_a=\{ z\in C : \varepsilon^a\leq |z|\leq \varepsilon_0\}$, where $a$
will be chosen as follows:

\medskip
\noindent {\bf Claim}.  For:
\begin{equation}\label{E:choice}
\frac{1+\kappa}{3-\kappa}> a > \frac{\kappa}{2-\kappa}\ ,
\end{equation}
there is a constant $B$ independent of $\ell$ so that:
$$
\int_{C_a} g^2 ds_0^2\leq \frac{B}{\varepsilon^{2\kappa}}\ .
$$

\bigskip
\noindent
Assuming this claim, we complete the proof.  By (\ref{E:L4bound}) and Prop.\ 
\ref{P:conicdegeneration} applied to (\ref{E:firstcoefficient}), we have:
$$
\int_{C_a}\frac{|c_{-1}|^4}{|z|^4}ds_0^2\leq \frac{B}{\varepsilon^{2\kappa}}
+\int_{C_a} \left(
\frac{|c_{-1}^{12}|^4|z|^{4\gamma}}{|z|^4}+\frac{|c_{-1}^{21}|^4|z|^{-4\gamma}}{|z|^4}          
\right) ds_0^2\ .
$$
Carrying out the integrals on both sides, this implies $|c_{-1}|^4\leq
B\left(\varepsilon^{(4-2\kappa)a-2\kappa}+\varepsilon^{4\gamma(1-a)}\right)$. 
 By the choice of
$a$ in (\ref{E:choice}),
 $(4-2\kappa)a-2\kappa>0$, so we are finished.

It remains to prove the claim.  First, by the expression for $|h^{-1}\partial_{A_0}h|^2$, note
that the terms in the series for $g$ involving the coefficients $c_m$ and $c_n$ are bounded
uniformly by $\varepsilon_0^{(|m|+|n|)b}$ for some $b>0$ and $|m|,|n|$ large, depending upon
$a$.   Therefore, to bound the integral of $g^2$, it suffices to bound the squares of the
individual terms.  Of the terms which appear, there are three types which need to be
estimated:
\begin{enumerate}
\item[I.]  $|c_m^{ij}|^2|c_n^{ij}|^2|z|^{2(m+n)+2(\gamma_i-\gamma_j)}$, where $m,n\leq -2$;
\item[II.]  $|c_{-1}^{21}|^2|c_n^{21}|^2|z|^{2(n-1)-2\gamma}$ and
$|c_{-1}^{12}|^2|c_n^{12}|^2|z|^{2(n-1)+2\gamma}$, where $n\leq -2$;
\item[III.] $|c_{-1}|^2|c_n^{ii}|^2|z|^{2(n-1)}$, where $n\leq -2$.
\end{enumerate}
Terms of type I give integrals of the form (using (i) above):
\begin{align*}
\int_{\varepsilon^a}^{\varepsilon_0}|c_m^{ij}|^2|c_n^{ij}|^2
r^{2(m+n)+2(\gamma_i-\gamma_j)-1+2\kappa} dr
&\leq
B|c_m^{ij}|^2|c_n^{ij}|^2\varepsilon^{(2(m+n)+2(\gamma_i-\gamma_j)+2\kappa)a}\\
&\leq B\varepsilon^{-2(m+n)(1-a)-2(\gamma_i-\gamma_j)(1-a)-2(2-\kappa a)}\ .
\end{align*}
But:
\begin{align*}
-2(m+n)(1-a)-2(\gamma_i-\gamma_j)(1-a)-2(2-\kappa a)&\geq
-2(m+n)(1-a)-2(1-a)-2(2-\kappa a)\\
&\geq -2(m+n+3)(1-a)\ ,
\end{align*}
so these terms are, in fact, bounded.  For type III, we apply the same
estimate to get:
$$
\int_{\varepsilon^a}^{\varepsilon_0}|c_{-1}|^2|c_n^{ii}|^2
r^{2(n-1)-1+2\kappa} dr
\leq
B |c_{-1}|^2 \varepsilon^{-2(n-1)(1-a)-2(2-\kappa a)}\ .
$$
Since $n\leq -2$, the exponent of $\varepsilon$ is $\geq -2\kappa$ by the assumption on $a$ in
(\ref{E:choice}).  Type II is similar to these two computations.
\end{proof}

Since the hypotheses of Prop.\ 
\ref{P:estimatetwo} are  satisfied by Cor.\  \ref{C:conicexpdecay} (2), and the uniform $L^2$ bound
on $|h^{-1/2}\partial_{A_0}h^{1/2}|$, this  completes the proof of (\ref{E:piconvergence}) and also
of the Main Theorem.

\end{document}